\input amstex
\documentstyle{amsppt}
\NoBlackBoxes
\nologo
 \magnification1200
\pageheight{9.2 true in}
\pagewidth{6.5 true in}

\topmatter
\title The two dimensional distribution of values of $\zeta(1+it)$
\endtitle
\author Youness Lamzouri
\endauthor
\address D{\'e}partment  de Math{\'e}matiques et Statistique,
Universit{\'e} de Montr{\'e}al, CP 6128 succ Centre-Ville,
Montr{\'e}al, QC  H3C 3J7, Canada
\endaddress
\email{Lamzouri{\@}dms.umontreal.ca}
\endemail
\thanks
\endthanks
\abstract We prove several results on the distribution
function of $\zeta(1+it)$ in the complex plane, that is the joint distribution function of $\arg\zeta(1+it)$ and $|\zeta(1+it)|$. Similar
results are also given for $L(1,\chi)$ (as $\chi$ varies over non-principal
characters modulo a large prime $q$).
\endabstract
\toc

\head  Introduction \page{1}
\endhead

\head  1. Detailed statement of results \page{4}
\endhead

\head 2. Approximations of $\zeta(1+it)$\page{7}
\endhead

\head 3. Estimates of sums of divisor functions\page{10}
\endhead

\head 4. Moments of $\zeta(1+it)$ \page{15}
\endhead

\head 5. Large values of $\zeta(1+it)$ in every direction\page{18}
\endhead

\head 6. Random Euler products and their distribution \page{24}
\endhead

\head 7. Fourier analysis on the $n$-dimensional torus
\page{28}
\endhead

\head 8. The normal distribution of $\arg\zeta(1+it)$\page{32}
\endhead

\head 9. Analogous results for $L(1,\chi)$ \page{35}
\endhead

\head  References  \page{40}
\endhead

\endtoc

\endtopmatter

\document
\head Introduction \endhead

\noindent The values of the Riemann zeta function and $L$-functions
at the edge of the critical strip $\text{Re} (s)=1$, have important
arithmetical consequences. The first one being the fact that

\noindent $\zeta(1+it)\neq 0$ implies the prime number theorem,
proved by Hadamard and de La Vall\'ee Poussin in 1896, that
$$ \pi(x)\sim \frac{x}{\log x}, \ \ \text{as} \ \ x\to\infty.$$
The second one is the class number formula, proved by Dirichlet in 1839,
which relates the class number of a quadratic extension of $\Bbb Q$
to the value of $L(1,\chi_d)$ where $d$ is the discriminant of the
field extension.

 The distribution of these values have been extensively
studied over the last decades. One can quote the work of
Granville-Soundararajan [10] in the case of $|\zeta(1+it)|$; Elliott ([7] and [8]), Montgomery-Vaughan [19] and
Granville-Soundararajan [11]
 in the case of Dirichlet $L$-functions
of quadratic characters $L(1,\chi_d)$; Duke [6] in the case of Artin $L$-functions, and the work of
Cogdell-Michel [4], Habsieger-Royer [12], Lau-Wu [14], Liu-Royer-Wu [17], Royer ([21] and [22]), and Royer-Wu ([23] and [24]) in the case of symmetric power
$L$-functions of $GL_2$-automorphic forms.

 We know that the Riemann zeta function $\zeta(s)$ has a
 conditionally convergent Euler product on Re$(s)=1$
$$\zeta(1+it)=\lim_{y\to\infty}\prod_{p\leq y}\left(1-\frac{1}{p^{1+it}}\right)^{-1}, \text{ if } t\gg 1.\tag {1}$$
\noindent In 1928, assuming the Riemann Hypothesis, Littlewood ([15]
and [16]) showed that one can truncate this product at $p\leq
\log^2t$ to obtain a good approximation for $\zeta(1+it)$, deducing
that $|\zeta(1+it)|\leq (2 e^{\gamma}+o(1)) \log_2 t$. (Throughout
$\log_j$ denotes the $j$-th iterated logarithm, so that $\log_1n
=\log n$ and $\log_j n=\log (\log_{j-1}n) $ for each $j\geq 2$). This shows that under the Riemann Hypothesis the sum $\sum_{p\geq y}1/p^{1+it}$ is small for $y\geq \log ^2 t.$ Moreover using Dirichlet's Theorem on diophantine approximation it is possible to make the sum $\sum_{p\leq \log t}1/p^{1+it}$ large, by choosing $t$ such that $p^{it}\approx 1$, for all the primes $p\leq \log t$. This enabled Littlewood ([15]
and [16])  to show the existence of arbitrarily
large $t$ for which $|\zeta(1+it)|\geq  (e^{\gamma}+o(1)) \log_2 t$. Furthermore it is widely believed that the sum $\sum_{\log t\leq p\leq \log ^2 t}1/p^{1+it}$ is small so that the truncated product up to
$\log t$ still serves as a good approximation for $\zeta(1+it)$:

 \proclaim {
Conjecture 1} As $t\to \infty$, we have
$$\zeta(1+it)\sim\prod_{p\leq \log
t}\left(1-\frac{1}{p^{1+it}}\right)^{-1}.
 $$
\endproclaim
\noindent One consequence of this conjecture is that $\max_{|t|\leq
T} |\zeta(1+it)| \sim e^{\gamma} \log_2 T.$ In
2003, Granville and Soundararajan [10] evaluate the frequency with
which such extreme values are attained, giving strong evidence for
the truth of Conjecture 1. More precisely if $$
  {\Phi}_T(\tau):= \frac 1T \text{meas} \{ t\in [T,2T]: \
\ |\zeta(1+it)| >e^{\gamma}\tau\},$$ then uniformly in the range
$1\ll \tau \le \log_2 T -\log_3 T$, they proved that
$$
{\Phi}_T(\tau) = \exp\left(-\frac{2e^{\tau-C-1}}{\tau} \left(1 +
O\left(\frac{1}{\tau^{\frac 12}}\right)\right)\right) ,\tag{2}
$$
where $$C=\int_0^2\log I_0(t)\frac{dt}{t^2} + \int_2^{\infty} (\log
 I_0(t)-t)\frac{dt}{t^2},\tag{3}$$
is a positive constant and $I_0(t): = \sum_{n=0}^{\infty}
(t/2)^{2n}/n!^2$ is the modified Bessel function of order $0$.

 The aim of this paper is to investigate the tail of the joint distribution
function of $|\zeta(1+it)|$ and $\arg\zeta(1+it)$ (where the latter is defined by continuous variation of the argument along the straight lines joining $2$, $2+it$ and $1+it$ starting with the value $0$):
 $$\Phi_T(\tau,\theta):=\frac{1}{T}\text{meas}\{t\in[T,2T]: |\zeta(1+it)|>e^{\gamma}\tau,
 \ |\arg\zeta(1+it)|>\theta\},$$
 for $\tau$ large and $\theta>0$ bounded. In the same range $1\ll\tau\leq \log_2
 T-\log_3 T$ as in (2), we show (in Theorem 1) that for any fixed
 $\theta>0$
 $$\Phi_T(\tau,\theta)=\exp\left(-e^{\tau(1+o_{\theta}(1))}\right),\tag{4}$$
so the proportion does not decay too fast. We can be more precise
showing (see Theorem 5 below), in the smaller range $1\ll\tau\leq
(\log_2 T)/2-2\log_3 T$, and $(\log \tau)\sqrt{\frac{\log_2
\tau}{\tau}}<\theta\ll 1$, that
$$
\Phi_T(\tau,\theta)=
\exp\left(-\frac{\displaystyle{e^{\tau+\frac{\theta^2\tau}{2\log\tau}
+O\left(\frac{\theta^2\tau}{\log^2\tau}\right)}}}{\tau}\right).\tag{5}$$
We do prove the implicit upper bound in the full range $1\ll\tau\leq
\log T-\log 10$ unconditionally, and that the lower bound holds in this
range assuming the Lang-Waldshmidt conjecture for linear forms in
logarithms (Conjecture 2 below). As a consequence of our result we
deduce that almost all values of $\zeta(1+it)$ with large norm are
concentrated near the positive real axis:

\proclaim {Corollary 1} As $\tau, T\to \infty$ with $\tau\leq \log_2T-\log_3 T$,  almost all values
$t\in [T,2T]$, with $|\zeta(1+it)|>e^{\gamma}\tau$, satisfy
$|\arg\zeta(1+it)|\leq (\log \tau)\sqrt{\log_2\tau/\tau}$. Moreover the set of exceptions has measure $\leq \exp(-\exp(\tau+(\log\tau\log_2\tau)/2))$.
\endproclaim
Also from the estimate (5), one can deduce that the larger the
arguments, the more it becomes rare to find values with large norm.
More precisely we have

\proclaim {Corollary 2} Let $\tau$, $\theta_1$ and $\theta_2$,
be in the range of validity of Theorem 5. If $\tau$ is large and

\noindent  $\theta_1>\theta_2 (1+c_5/\log\tau)$, where $c_5$ is a
suitably large constant, then
$$ \Phi_T(\tau,\theta_1)=o(\Phi_T(\tau,\theta_2)), \text{ as } \tau, T\to \infty.$$

\endproclaim

Let $\tau\leq \log_2 T$ be a large real number. Another interesting question is to understand  the behavior of the argument of $\zeta(1+it)$ for  $t$ with $|\zeta(1+it)|\approx e^{\gamma} \tau$.
The appearance of  the factor $(\theta^2/2)\tau/\log \tau$
in (5), may suggests a normal behavior in the argument $\theta$. Indeed
 we evaluate the characteristic function of
$\arg\zeta(1+it)$ with an appropriate weight, and show that these
arguments should be distributed according to a normal law of mean
$0$ and variance $\log(\tau-1-C)/2e^{\tau-1-C}$ (see Theorem 6
below).

We will introduce a  random model for the values $\zeta(1+it)$ : Let
$\{X(p)\}_{p \ \text{prime}}$ be a set of independent random
variables, uniformly distributed on the unit circle ${\Bbb U}$, and
define the ``random Euler product''
$$L(1,X)=\lim_{y\to\infty}\prod_{p\leq
y}\left(1-\frac{X(p)}{p}\right)^{-1}, \ (\text{these products
converge with probability} \ 1).
$$
Our strategy is to compare the distribution of the values of
$\zeta(1+it)$ with the distribution of $L(1,X)$. For example we show
in Theorem 2 below, that large complex moments of $\zeta(1+it)$ and
$L(1,X)$ are roughly equal (Granville and Soundararajan
(unpublished) proved an analogous result for $L(1,\chi)$, see
Theorem B in section 9). Therefore we study this probabilistic model
closely (Theorem 3) and deduce results on the distribution of
$\zeta(1+it)$ (Theorem 5).

The results proved here carry over to $L(1,\chi)$ (where $\chi$
varies over non-principal characters modulo a large prime $q$)
without any difficulty. We discuss these results in section 9.

\noindent { \bf Acknowledgments.} I sincerely thank my advisor, Professor Andrew Granville,
 for suggesting this problem and for all his advice and encouragement. I would also thank Professor K. Soundararajan for valuable discussions.

\head 1. Detailed statement of results \endhead

\noindent First we define
 $$ \zeta(1+it,y):=\prod_{p\leq
 y}\left(1-\frac{1}{p^{1+it}}\right)^{-1}, \ \ \text{and} \ \ R_y:=\prod_{p\leq
 y}\left(1-\frac{1}{p}\right)^{-1}.$$
To exhibit large values of $\zeta(1+it)$ in any given direction
 $\arg z=\theta$, we  first approximate $\zeta(1+it)$ by short Euler products
 $\zeta(1+it,y)$
 (which
 is possible for almost all $t\in[T,2T]$ by Lemma 2.4 below in the range $1\ll y\leq \log T$), then we try to find
 many values $t\in[T,2T]$ for which
 $$ \zeta(1+it,y)\approx
 e^{i\theta}R_y.$$
 To do so we use a biased method of moments, which we describe below. The first step is to note that the following inequality

$$ \left|\zeta(1+it,y)+ e^{i\theta}R_y\right|\geq (2-\epsilon) R_y,\tag{1.1}$$
implies
$$ \zeta(1+it,y)=e^{i\theta}R_y\left(1+O\left(\sqrt{\epsilon}\right)\right).$$
 This follows from the fact that  $ \left|\zeta(1+it,y)\right|\leq
 \left|e^{i\theta}R_y\right|,$ and noting that for a complex number $|z|\leq 1$, with
 $|z+1|\geq 2-\epsilon$, one can easily show that
$z=1+O(\sqrt{\epsilon})$. To prove (1.1) we can try to have a good lower bound for the moments
 $$
 \align
 \frac{1}{T}\int_T^{2T} &\left|\zeta(1+it,y)+
 e^{i\theta}R_y\right|^{2k}dt\\
 &=\sum_{0\leq l,m\leq k}\binom kl\binom km
R_y^{2k-l-m}e^{i\theta(m-l)}\frac{1}{T}\int_T^{2T}
\zeta(1+it,y)^l\zeta(1-it,y)^mdt.\tag{1.2}\\
\endalign
$$
In general, we can estimate these moments if the central terms $m=l$ constitute
the main term, (since for most cases it's difficult to handle the
non-central ones). However this is not the case here. In fact if
$y\leq (\log T)^2$, and $m,l\leq \log T/(25\log_2 T\log_3 T)$, then
by Theorem 4.1 below, we have
$$\frac{1}{T}\int_T^{2T}
\zeta(1+it,y)^l\zeta(1-it,y)^mdt = \sum_{n\in
S(y)}\frac{d_l(n)d_m(n)}{n^2} + o(1),\tag{1.3}$$ and so by
Proposition 3.2 below one can see that some non-central terms have the
same order as the central ones. Therefore it seems difficult to
estimate (1.2), because of the oscillation of $e^{i\theta(m-l)}.$
To handle this, we slightly modify the moments. Indeed,
instead of working with $\zeta(1+it,y)$, we search for some
completely multiplicative function $f(n)$ with values on the unit
circle $\Bbb U$, such that
$$ \prod_{p\leq
 y}\left(1-\frac{f(p)}{p^{1+it}}\right)^{-1}=R_y(1+O(\epsilon))\Longleftrightarrow
 \zeta(1+it)= e^{i\theta}R_y(1+O(\epsilon)).$$
In this case the non-central terms in
$$ \frac{1}{T}\int_T^{2T} \left|\prod_{p\leq y}\left(1-\frac{f(p)}{p^{1+it}}\right)^{-1}+
 R_y\right|^{2k}dt,\tag{1.4}$$ will be positive (by Theorem 4.1), and the
 central ones will give the lower bound we search for. In fact it
 turns out that the function we need, satisfies $f(p)=e^{-i\psi}$
 for all $p\leq y$, where $\psi=\theta/\log_2 y$.

 Using this method
we can prove the existence of large values of
$\zeta(1+it)$ in each given direction $\arg z=\theta.$ Indeed we
prove

\proclaim {Theorem 1} Let $T$ be large, and fix $\theta\in
(-\pi,\pi]$. If $1\ll y\leq \log T/\log_2 T$ is a real number, let
$M(\theta,y)$ be the measure of values $t\in[T,2T]$, for which
$$ \zeta(1+it)=e^{i\theta}\prod_{p\leq
y}\left(1-\frac{1}{p}\right)^{-1}\left(1+O\left(\frac{1}{\log_2
y}\right)\right).$$ Then there exist two positive constants
$c_1,c_2$ (depending on the constant in the $O$) for which
$$T\exp\left(-y^{1-c_2/(\log_2
y)^2}\right)\leq M(\theta,y)\leq T\exp\left(-y^{1-c_1/\log_2
y}\right).$$
\endproclaim
Granville and Soundararajan (unpublished) used a different method to
prove the existence of large values (and small values) of
$L(1,\chi)$ in every direction (see Theorem A of section 9). However
they only got a lower bound for the measure, and their bound is less
strong than what we obtain in Theorem 1.

Let $z$ be a complex number. We define the ``$z$th divisor function'' $d_z(n)$, to be the multiplicative function such that $d_z(p^a)=\Gamma(z+a)/\Gamma(z)a!$, for any prime $p$ and any integer $a\geq 0$. Then $d_z(n)$ is the coefficient of the Dirichlet series $\zeta(s)^z$ for Re$(s)>1$. Therefore for the random variables $\{X(p)\}_{p \text{ prime }}$ we have (with probability $1$) that

$$ L(1,X)^z=\sum_{n=1}^{\infty}\frac{d_z(n)X(n)}{n},$$
where $X(n)=\prod_{i=1}^k X(p_i)^{a_i},$ if $n=\prod_{i=1}^k p_i^{a_i}$. If $Y$ is a random variable on a probability space $(\Omega, \mu)$ we define its expectation by  ${\Bbb E}(Y)=\int_{\Omega}Yd\mu.$
 Therefore
${\Bbb E}(X(n)\overline{X(m)})=1$  if $n=m$ and vanishes otherwise.
Thus for any complex numbers
$z_1$ and $z_2$, we have
$${\Bbb
E}\left(L(1,X)^{z_1}\overline{L(1,X)}^{z_2}\right)=\sum_{n=1}^{\infty}\frac{d_{z_1}(n)d_{z_2}(n)}{n^2}.$$

The idea of using a probabilistic model appears previously in the work of
Montgomery-Vaughan [19], Granville-Soundararajan [11], and
Cogdell-Michel [4]. Indeed in each of these cases an adequate
probabilistic model was constructed to understand the distribution
of appropriate $L$-functions. To convince ourselves that it is the right model to use, we
evaluate high complex moments of $\zeta(1+it)$ and found that
\proclaim { Theorem 2} Uniformly for all complex numbers $z_1,z_2$
in the region

\noindent $|z_1|,|z_2|\leq \displaystyle{\log T/(50(\log_2 T)^2)}$,
we have
$$ \frac{1}{T}\int_T^{2T}\zeta(1+it)^{z_1}\zeta(1-it)^{z_2}dt ={\Bbb
E}\left(L(1,X)^{z_1}\overline{L(1,X)}^{z_2}\right)
+ O\left(\exp\left(-\frac{\log T}{2\log_2 T}\right)\right),$$
\endproclaim
 \noindent  Using a combinatorial argument, Granville and Soundararajan [10],
get a better  result (in the uniformity of the range of moments) in
the special case where $z_2=z_1=k\in{\Bbb Z}$.

This result motivated us to study the distribution of the random
Euler products $L(1,X)$. For $\tau, \theta>0$, define
$$ \Phi(\tau,\theta):= \text{Prob} (|L(1,X)|>e^{\gamma}\tau,
 \ |\arg L(1,X)|>\theta).$$
 A close study of this model allowed us to find a precise estimate
 for this distribution function. Indeed we prove the following

\proclaim{ Theorem 3} For $\tau>0$ large and
$(\log \tau)\sqrt{\frac{\log_2 \tau}{\tau}}
<\theta\ll 1$, we have
$$
\Phi(\tau,\theta)=\exp\left(-\frac{\displaystyle{e^{\tau+\frac{\theta^2\tau}{2\log\tau}
+O\left(\frac{\theta^2\tau}{\log^2\tau}\right)}}}{\tau}\right).$$
\endproclaim

Let $p_j$ denotes the $j$-th smallest prime number. To prove an analogous formula for $\Phi_T(\tau,\theta)$, we studied
the behavior of the vector $V(t):=(p_1^{it},p_2^{it},...,p_N^{it})$, in
the torus ${\Bbb T}^N:=({\Bbb R}/{\Bbb Z})^N$, for $t\in[T,2T]$, as $T\to \infty$. In fact we believe that these values should be equidistributed on ${\Bbb T}^N$ for $N=\pi(y)$ and  $y\leq  (1+o(1))\log T$. In [1],
Barton, Montgomery and Vaaler an hold is constructed trigonometric polynomials
in $N$ variables, which give a sharp approximation to the
characteristic function of a cartesian product of $N$ open intervals
(see Theorem 7.1). These polynomials are the analogue of Selberg
polynomials in $1$ variable (see [18]). Using this construction and Fourier analysis on ${\Bbb T}^N$,
 we show in Theorems 4A and 4B below, that these values are equidistributed on ${\Bbb T}^N$, for $y\leq \sqrt{\log T}/(\log_2 T)^2$
 unconditionally, and for $y\leq (\log T)/10$ under a conjecture on linear forms in logarithms,
 formulated by Lang and Waldshmidt [13, Introduction to chapter X and XI, p. 212]:

\proclaim{Conjecture 2} Let $b_i$ be integers, and $a_i$ be positive
integers for which $\log a_i$
 are linearly independent over ${\Bbb Q}$. We let $B_j=\max\{|b_j|,1\}$, and $B=\max_{1\leq j\leq n} B_j$. Then for any $\epsilon>0$,
  there exists a positive
constant $c(\epsilon)$, such that
$$ |b_1\log a_1 + b_2\log a_2+...+b_n\log
a_n|>\frac{c(\epsilon)^nB}{(B_1...B_na_1...a_n)^{1+\epsilon}}.
$$
\endproclaim
More precisely we prove

\proclaim {Theorem 4A} Let $2<y$ be a real number. For each $1\leq
j\leq \pi(y)$, let $I_j\subset (0,1)$ be an open interval of length
$\delta_j>0$. Define
$$ M(I_1,...,I_{\pi(y)})=M:=\text{meas}\left\{ t\in [T,2T] : \left\{ \frac{t\log p_j}{2\pi} \right\}\in
I_j, \text{ for all } 1\leq j\leq \pi(y) \right\},$$ where $\{\cdot\}$ denotes the fractional part. Then
$$M\sim T\prod_{j\leq \pi(y)}\delta_j,$$
uniformly for $y\leq \sqrt{\log T}/(\log_2 T)^2$, and
$\delta_j>(\log_2 T)^{-5/3}$.
\endproclaim
\proclaim {Theorem 4B} Assume Conjecture 2. Then with the same
notations as Theorem 4A, we have
$$M\sim T\prod_{j\leq \pi(y)}\delta_j,$$
uniformly for $y\leq (\log T)/10$, and $\delta_j>(\log T)^{-3/2}$.
\endproclaim

Following the proof of Theorem 3 and using Theorems 4A and 4B we deduce

\proclaim{ Theorem 5} Let $T>0$ be large. There exists two positive
constants $c_3$ and $c_4$ for which
$$
\Phi_T(\tau,\theta)\leq
\exp\left(-\frac{\displaystyle{e^{\tau+\frac{\theta^2\tau}{2\log\tau}
-c_3\frac{\theta^2\tau}{\log^2\tau}}}}{\tau}\right),$$
uniformly for $1\ll\tau\leq \log_2 T$, and $(\log \tau)\sqrt{\frac{\log_2
\tau}{\tau}}<\theta\ll 1$. And
$$
\Phi_T(\tau,\theta)\geq
\exp\left(-\frac{\displaystyle{e^{\tau+\frac{\theta^2\tau}{2\log\tau}
+c_4\frac{\theta^2\tau}{\log^2\tau}}}}{\tau}\right),$$
uniformly for $(\log \tau)\sqrt{\frac{\log_2
\tau}{\tau}}<\theta\ll 1$, and $1\ll\tau\leq (\log_2 T)/2-2\log_3 T$
unconditionally, and for
 $1\ll\tau\leq \log_2 T-\log 10$ if we assume Conjecture 2.
\endproclaim

We now turn our attention to the behavior of $\arg\zeta(1+it)$ when the norm is large,
 that is when $|\zeta(1+it)|\approx e^{\gamma}\tau$ with $\tau\leq (1+o(1))\log_2 T$.
 We compute the characteristic function of these arguments with a natural weight,
  and use the Berry-Esseen Theorem ([2], [9]) to prove the following

\proclaim{ Theorem 6 } Let $T>0$ be large,  $1\ll\tau\leq
\log_2T-3\log_3T$ a real number, $\epsilon=\tau^{-1/5}$ and
$k=e^{\tau-1-C}$, where $C$ is defined by (3). Let
$$ \Omega_T(\tau):=\{t\in [T,2T]: \ e^{\gamma}(\tau-\epsilon)\leq
|\zeta(1+it)|\leq e^{\gamma}(\tau +\epsilon)\},$$ and for a real
number $x$, let
$$ \Lambda_T(\tau,x):=\left\{t\in\Omega_T(\tau): \
\displaystyle{\frac{\arg\zeta(1+it)}{\sqrt{\frac{\log(\tau-1-C)}{2e^{\tau-1-C}}}}}<x\right\},
 \text{ and } \
\nu_{T,\tau}(x):=\frac{\displaystyle{\int_{\Lambda_T(\tau,x)}|\zeta(1+it)|^{2k}dt}}
{\displaystyle{\int_{\Omega_T(\tau)}|\zeta(1+it)|^{2k}dt}}.$$ Then
we have
$$\nu_{T,\tau}(x)=\frac{1}{\sqrt{2\pi}}\int_{-\infty}^{x}e^{-y^2/2}dy+O_x\left(\frac{1}{\sqrt{\log\tau}}\right).$$
\endproclaim

\head 2. Approximations of $\zeta(1+it)$\endhead

\subhead 2.1 Short Euler product approximation \endsubhead

\smallskip

\noindent In this section we approximate $\zeta(1+it)$ by a short
Euler product of length $y\leq \log T$, for almost all $t\in
[T,2T]$. The main idea is to show that this is possible if
$\zeta(s)$ has no zeros far from the critical line, then to use a
classical zero-density estimate (there are  few such zeros) to see
that we can almost surely avoid there zeros. The material of this
section is classical, and it's essentially proved in Granville and Soundararajan [10] (see sections 2 and 3).

 \proclaim{Lemma 2.1  ([10, Lemma 1])} Let $y\ge 2$ and $|t| \ge y+3$ be real numbers.
Let $\frac{1}{2} \leq \sigma_0 <1$ and suppose that the rectangle
$\{ s: \ \ \sigma_0 <\text{Re}(s) \leq 1, \ \ |\text{Im}(s) -t| \leq
y+2\}$ does not contain any zeros of $\zeta(s)$. Then if  $\sigma_0 <
\sigma \le 2$ and $|x-t|\leq y$ we have
$$
|\log \zeta(\sigma +ix)| \ll \log |t|  \log (e/(\sigma -\sigma_0)).
$$
Moreover, if $\sigma_0 < \sigma \leq 1$ then
$$
\log \zeta(\sigma+it)= \sum_{n=2}^{y}
\frac{\Lambda(n)}{n^{\sigma+it} \log n} + O\left( \frac{\log
|t|}{(\sigma_1-\sigma_0)^2}y^{\sigma_1-\sigma}\right),
$$
where  $\sigma_1 = \min(\sigma_0+\frac{1}{\log y},
\frac{\sigma+\sigma_0}{2})$.
\endproclaim

\noindent From this result, we deduce

\proclaim{Lemma 2.2 ([10, Lemma 2]) } Let $\frac{1}{2} < \sigma \le 1$ be fixed, $T$
large and $3<y<T/2$ be a real number. We have
$$
\log \zeta(\sigma+it) = \sum_{n=2}^{y}
\frac{\Lambda(n)}{n^{\sigma+it} \log n} + O( y^{(\frac 12-\sigma)/2}
\log^3 T)
$$
for all $t\in (T,2T)$ except for a set of measure $\ll
T^{5/4-\sigma/2} y (\log T)^5$.
\endproclaim

\demo{Proof}  This follows from combining the classical zero-density
estimate

\noindent $N(\sigma_0,T) \ll T^{3/2-\sigma_0} (\log T)^5$ (see
Theorem 9.19 A of [25]) and  Lemma 2.1 (taking $\sigma_0 = (1/2
+\sigma)/2$ there).
\enddemo
\noindent To obtain an approximation by shorter Euler products, we
need a large sieve type inequality

\proclaim {Lemma 2.3 ([10, Lemma 3])}  Let $2\le y\le z$ be real numbers. For
arbitrary complex numbers $x(p)$ we have
$$
\frac{1}{T} \int_{T}^{2T} \left| \sum_{y\le p\le z}
\frac{x(p)}{p^{it}} \right|^{2k}dt \ll \left(k\sum_{y\le p \le z}
|x(p)|^2 \right)^k + T^{-\frac{2}{3}} \left( \sum_{y\le p\le z}
|x(p)|\right)^{2k}
$$
uniformly for all integers $1\le k \le \log T/(3\log z)$.
\endproclaim

 We define $\zeta(s,y):=\prod_{p\leq
y}\left(1-p^{-s}\right)^{-1}$, and using the Lemmas above, we prove
the following key Lemma

\proclaim {Lemma 2.4} Let $T>0$ be a large real number, and
$A(t)\leq \log t$ be a slowly increasing function which tends to
$\infty$  with $t$. Then, uniformly for $y\leq \log T$, we have
$$ \zeta(1+it)=\zeta(1+it,y)\left(1+O\left(\frac{1}{A(y)}\right)\right),$$
for all $t\in [T,2T]$ except a set of measure
$$\ll T\exp\left(-\log\left(\frac{300\log^2 y}{A(y)^2}\right)\frac{y}{300\log y}\right).
$$
\endproclaim
\demo{Proof} Let $z=(\log T)^{100}$, we deduce from Lemma 2.2 that
$$ \zeta(1+it)=\zeta(1+it,z)\left(1+O\left(\frac{1}{\log
T}\right)\right),\tag{2.1}$$ for all $t\in [T,2T]$ except a set of
measure at most $T^{4/5}$. Applying Lemma 2.3 with $x(p)=1/p$, we
get
$$ \frac{1}{T}\int_T^{2T}\left|\sum_{y\leq p\leq
z}\frac{1}{p^{1+it}}\right|^{2k}dt\ll \left(k\sum_{y\leq p\leq
z}\frac{1}{p^{2}}\right)^k+ T^{-2/3}\left(\sum_{y\leq p\leq
z}\frac{1}{p}\right)^{2k},$$ for any integer $1\leq k\leq \log
T/3\log z$. We choose $k=[y/(300\log y)]$, which implies that
$$ \frac{1}{T}\int_T^{2T}\left|\sum_{y\leq p\leq
z}\frac{1}{p^{1+it}}\right|^{2k}dt\ll\left(\frac{1}{300 \log^2
y}\right)^k.\tag{2.2}$$ Let
$M=\text{meas}\{t\in[T,2T]:\displaystyle{\left|\sum_{y\leq p\leq
z}\frac{1}{p^{1+it}}\right|>\frac{1}{A(y)}}\}.$ From (2.2) we get
$$ \frac{M}{T}\left(\frac{1}{A(y)}\right)^{2k}\ll \left(\frac{1}{300
\log^2 y}\right)^k,$$ which implies that  $$M\ll
T\exp\left(-\log\left(\frac{300\log^2
y}{A(y)^2}\right)\frac{y}{300\log y}\right).$$ To complete the proof
one may check that for all $t\in[T,2T]$, except a set of measure
$M$, we have
$$
\align \zeta(1+it,z)&=\zeta(1+it,y)\exp\left(-\sum_{y\leq p\leq
z}\left(\frac{1}{p^{1+it}}+O\left(\frac{1}{p^2}\right)\right)\right)\\
&=\zeta(1+it,y)\left(1+O\left(\frac{1}{A(y)}\right)\right).\\
\endalign
$$
\enddemo

\subhead 2.2 Smooth Dirichlet series approximation of
$\zeta(1+it)^z$\endsubhead

\smallskip

\noindent To prove Theorem 2, we need the following Lemma, which
corresponds to Lemma 2.3 of Granville-Soundararajan [11].

\proclaim {Lemma 2.5} Let $t$ be large, and $z$ be any complex number
with $|z|\leq \log^2 t$. Define $Z=\exp((\log t)^{10})$. Then
$$
\zeta(1+it)^z=\sum_{n=1}^{\infty}\frac{d_z(n)}{n^{1+it}}e^{-n/Z}+O\left(\frac{1}{t}\right).$$
\endproclaim

\demo{Proof} Since $\frac{1}{2\pi
i}\int_{1-i\infty}^{1+i\infty}y^s\Gamma(s)ds=e^{-1/y}$, we have
$$ \frac{1}{2\pi i} \int_{1-i\infty}^{1+i\infty}
\zeta(1+it+s)^zZ^s\Gamma(s) ds
=\sum_{n=1}^{\infty}\frac{d_z(n)}{n^{1+it}}e^{-n/Z}.\tag{2.3}$$ We
shift the line of integration to the contour $s=-C(x)+ix$ where
$C(x):=c/(2\log(|x|+2))$, and $c>0$ is chosen so that $\zeta(s)$
have no zeros in the region where we shift the contour (this is
possible by the classical theorem on the zero free region of
$\zeta(s)$). We encounter a pole at $s=0$, which leaves the residue
$\zeta(1+it)^z$. Applying Lemma 2.1 with $\sigma_0=1-4C(x)/3$ and
$y=2$ gives $|\log \zeta(s)|\ll 1/C(x)^2$ and so the integral along
the new contour is
$$\ll
\int_{-\infty}^{\infty}Z^{-C(x)}e^{O(|z|/C(x)^2)}|\Gamma(-C(x)+ix)|dx\ll\frac{1}{t},$$
by Stirling's formula. This completes the proof.

\enddemo

\head 3. Estimates of sums of divisor functions  \endhead

\noindent In this section we prove two results on sums of the
divisor function $d_k(n)$. The advantage of our results is the
uniformity on $k$. We begin by proving the following proposition on
the estimates of such sums in short intervals, which we shall use
later in the proof of Theorem 2

\proclaim{Proposition 3.1} Let $T>0$ be a large real number, and
$k\leq \log T/4(\log_2 T)^2$ a positive integer. Define $Z=\exp((\log
T)^{10})$ and $y=\exp(\log T/\log_2 T)$. Then
$$ \sup_{m>\sqrt{T}}\sum_{m<n<m+\frac{m}{\sqrt{T}}}\frac{d_k(n)}{n}e^{-n/Z}
\leq \frac{(\log 3Z)^{k}}{y}.$$
\endproclaim

\demo{ Proof} We prove this by induction on $k$. If $k=1$ then
$$\sup_{m>\sqrt{T}}\sum_{m<n<m+\frac{m}{\sqrt{T}}}\frac{e^{-n/Z}}{n}\leq
\sup_{m>\sqrt{T}}\frac{1}{m}\left(\frac{m}{\sqrt{T}}+1\right)\leq
\frac{2}{\sqrt{T}}\leq \frac{\log 3Z}{y}.$$ Now suppose the result
true for $k-1$. K. Norton [20] proved that
$$ \log d_k(n)\leq \frac{\log n\log k}{\log\log n}\left(1+\frac{\log
\log\log n}{\log\log n}\left(1+O\left(\frac{1}{\log\log
n}\right)\right)\right),$$ uniformly for $k\leq \log n/(\log\log
n)^2$, if $n$ is large enough. Thus
$$
\align
&\sup_{\sqrt{T}<m<y\sqrt{T}}\sum_{m<n<m+\frac{m}{\sqrt{T}}}\frac{d_k(n)}{n}e^{-n/Z}\leq
\sup_{\sqrt{T}<m<y\sqrt{T}}\sum_{m<n<m+y} \frac{d_k(n)}{n}\\
&\leq\sup_{\sqrt{T}<m<y\sqrt{T}}\left(y\max_{m<n<m+y}\frac{d_k(n)}{n}\right)\leq \frac{1}{y}.\tag{3.1}\\
\endalign
$$
Now for $m>y\sqrt{T}$, we have
$$ \sum_{m<n<m+\frac{m}{\sqrt{T}}}\frac{d_k(n)}{n}e^{-n/Z}=
\sum_{m<n<m+\frac{m}{\sqrt{T}}}\frac{e^{-n/Z}}{n}\sum_{dr=n}d_{k-1}(r).$$
We divide the above sum into two parts: $S_1$ for $d>y$ (which
implies that $r\leq \frac{2m}{y})$, and $S_2$, for $d\leq y$. We
have then
$$ S_1\leq\sum_{r\leq \frac{2m}{y}}
\frac{d_{k-1}(r)}{r}e^{-r/Z}\sum_{\frac{m}{r}<d<\frac{m}{r}+\frac{m}{r\sqrt{T}}}\frac{1}{d}\leq
\sum_{r\leq \frac{2m}{y}}
\frac{d_{k-1}(r)}{r}e^{-r/Z}\left(\frac{1}{\sqrt{T}}+\frac{2}{y}\right).
$$
For $j\in {\Bbb N}$, we have that  $d_j(n)e^{-n/Z}\leq
e^{j/Z}\sum_{a_1...a_j=n}e^{-(a_1+...+a_j)/Z}$, and so
$$ \sum_{n=1}^{\infty}\frac{d_j(n)}{n}e^{-n/Z}\leq \left(e^{1/Z}
\sum_{a=1}^{\infty}\frac{e^{-a/Z}}{a}\right)^j\leq (\log
3Z)^j.\tag{3.2}$$ This implies
$$ S_1\leq \frac{3(\log 3Z)^{k-1}}{y}\leq \frac{(\log 3Z)^{k}}{3 y}.\tag{3.3}$$
Moreover
$$ S_2\leq \sum_{d\leq
y}\frac{1}{d}\sum_{\frac{m}{d}<r<\frac{m}{d}+\frac{m}{d\sqrt{T}}}\frac{d_{k-1}(r)}{r}e^{-r/Z}.$$
Since $m>y\sqrt{T}$, and $d\leq y$, we get $m/d>\sqrt{T}$. By our
induction hypothesis we deduce that
$$
\sum_{\frac{m}{d}<r<\frac{m}{d}+\frac{m}{d\sqrt{T}}}\frac{d_{k-1}(r)}{r}e^{-r/Z}
\leq\sup_{s>\sqrt{T}}\sum_{s<r<s+\frac{s}{\sqrt{T}}}\frac{d_{k-1}(r)}{r}e^{-r/Z}\leq
\frac{(\log 3Z)^{k-1}}{y}.$$ Finally we have
$$ S_2\leq \frac{(\log 3Z)^{k-1}\log y}{y}\leq \frac{(\log
3Z)^{k}}{3y}.\tag{3.4}$$ Now combining (3.1), (3.3) and (3.4) gives
the result.
\enddemo

 The key ingredient of the proof of Theorem 6, is to
understand the ratio
$$\sum_{n=1}^{\infty}\frac{d_{k-r}(n)d_{k+r}(n)}{n^2}\Big/\sum_{n=1}^{\infty}\frac{d_k(n)^2}{n^2},$$
for large $k$ and $r$. To this end we prove the following result

\proclaim{ Proposition 3.2} Let $k$ be a large real number and

\noindent $c_0=\lim_{x\to\infty}\left(\sum_{p\leq x}1/p-\log_2
x\right)$. Then uniformly for $|r|\leq \sqrt{k}$, we have

$$\sum_{n\geq 1}\frac{d_{k-r}(n)d_{k+r}(n)}{n^2}=
\exp\left(-r^2\frac{\log_2 k}{k}-\frac{c_0r^2}{k}
+O\left(\frac{r^2}{k\sqrt{\log
k}}+\frac{r^4}{k^2}\right)\right)\sum_{n\geq
1}\frac{d_{k}^2(n)}{n^2}.$$
\endproclaim

\noindent First, we remark that
$$
\align
&\frac{1}{2\pi}\int_{-\pi}^{\pi}\left|1-\frac{e^{i\theta}}{p}\right|^{-2k}d\theta=
\frac{1}{2\pi}\int_{-\pi}^{\pi}\left(1-\frac{e^{i\theta}}{p}\right)^{-k}\left(1-\frac{e^{-i\theta}}{p}\right)^{-k}d\theta\\
&=\frac{1}{2\pi}\int_{-\pi}^{\pi}\sum_{a=0}^{\infty}\frac{d_k(p^a)e^{ia\theta}}{p^a}
\sum_{b=0}^{\infty}\frac{d_k(p^b)e^{-ib\theta}}{p^b}d\theta=\sum_{a=0}^{\infty}\frac{d_k^2(p^a)}{p^{2a}}.\tag{3.5}\\
\endalign
$$
Analogously we have
$$
\align &\sum_{a=0}^{\infty}\frac{d_{k-r}(p^a)d_{k+r}(p^a)}{p^{2a}}=
\frac{1}{2\pi}\int_{-\pi}^{\pi}\left|1-\frac{e^{i\theta}}{p}\right|^{-2k}\left(1-\frac{e^{i\theta}}{p}\right)^r
\left(1-\frac{e^{-i\theta}}{p}\right)^{-r} d\theta\\
&=\frac{1}{4\pi}\int_{-\pi}^{\pi}\left|1-\frac{e^{i\theta}}{p}\right|^{-2k}
\left(\left(1-\frac{e^{i\theta}}{p}\right)^r
\left(1-\frac{e^{-i\theta}}{p}\right)^{-r}
+\left(1-\frac{e^{-i\theta}}{p}\right)^r
\left(1-\frac{e^{i\theta}}{p}\right)^{-r} \right)d\theta\\
&=\frac{1}{2\pi}\int_{-\pi}^{\pi}\left|1-\frac{e^{i\theta}}{p}\right|^{-2k}\cos\left(
2r\arg\left(1-\frac{e^{i\theta}}{p}\right)\right)d\theta,\tag{3.6}\\
\endalign$$
since the series
$\sum_{a=0}^{\infty}d_{k-r}(p^a)d_{k+r}(p^a)/p^{2a}$ is real. The
proof will rely on these two identities. The last ingredient we need
is the following Lemma
\proclaim{ Lemma 3.3} Let $k>0$ be a large
real number. Suppose that $p=o(k)$ as $k\to\infty$, and let
$\epsilon= 4\sqrt{\frac{p}{k}\log\left(\frac{k}{p}\right)}$. Then we
have
$$
\frac{1}{2\pi}\int_{-\epsilon}^{\epsilon}\left|1-\frac{e^{i\theta}}{p}\right|^{-2k}d\theta
= \left(\frac{1}{2\pi}\int_{-\pi}^{\pi}\left|1-\frac{e^{i\theta}}{p}
\right|^{-2k}d\theta\right)\left(1+O\left(\frac{p^4}{k^4}\right)\right).$$

\endproclaim

\demo{Proof} First we observe that
$$
\align \frac{1}{2\pi}\int_{-\pi}^{\pi}\left|1-\frac{e^{i\theta}}{p}
\right|^{-2k}d\theta &=
\frac{1}{2\pi}\int_{-\pi}^{\pi}\left(1-\frac{2\cos\theta}{p}+\frac{1}{p^2}\right)^{-k}d\theta\\
&=
\frac{1}{\pi}\int_{0}^{\pi}\left(1-\frac{2\cos\theta}{p}+\frac{1}{p^2}\right)^{-k}d\theta.\\
\endalign
$$
Now
$$
\align
\frac{1}{\pi}\int_{\epsilon}^{\pi/2}\left(1-\frac{2\cos\theta}{p}+\frac{1}{p^2}\right)^{-k}d\theta
&\leq
\frac{1}{\sin\epsilon}\left(\frac{1}{\pi}\int_{\epsilon}^{\pi/2}
\sin\theta\left(1-\frac{2\cos\theta}{p}+\frac{1}{p^2}\right)^{-k}d\theta\right)\\
&=\frac{1}{\sin\epsilon}\left(-\frac{p}{2\pi(k-1)}
\left(1-\frac{2\cos\theta}{p}+\frac{1}{p^2}\right)^{-k+1}\Bigg|_{\epsilon}^{\pi/2}\right)\\
&\leq\frac{p}{\pi\epsilon(k-1)}\left(1-\frac{1}{p}\right)^{-2k+2}\exp\left(-\frac{\epsilon^2k}{p}
+O\left(\frac{\epsilon^4k}{p}\right)\right).\\
\endalign
$$
Moreover since
$$
\frac{1}{\pi}\int_{\pi/2}^{\pi}\left(1-\frac{2\cos\theta}{p}+\frac{1}{p^2}\right)^{-k}d\theta
\leq \frac{1}{2}\left(1+\frac{1}{p^2}\right)^{-k},$$ then
$$E:=\frac{1}{\pi}\int_{\epsilon}^{\pi}\left(1-\frac{2\cos\theta}{p}+\frac{1}{p^2}\right)^{-k}d\theta
\leq
\frac{2p}{\pi\epsilon(k-1)}\left(1-\frac{1}{p}\right)^{-2k+2}\exp\left(-\frac{\epsilon^2k}{p}
+O\left(\frac{\epsilon^4k}{p}\right)\right).$$  Let $0<
\delta<\epsilon$ be a small real number, to be chosen later.
We have
$$
\align
I:=\frac{1}{\pi}\int_{0}^{\epsilon}\left(1-\frac{2\cos\theta}{p}+\frac{1}{p^2}\right)^{-k}d\theta
&\geq
\frac{1}{\pi}\int_{0}^{\delta}\left(1-\frac{2\cos\theta}{p}+\frac{1}{p^2}\right)^{-k}d\theta\\
&\geq\frac{1}{\sin\delta}\left(-\frac{p}{2\pi(k-1)}
\left(1-\frac{2\cos\theta}{p}+\frac{1}{p^2}\right)^{-k+1}\Bigg|_{0}^{\delta}\right)\\
&\geq
\frac{p}{2\pi\delta(k-1)}\left(1-\frac{1}{p}\right)^{-2k+2}\left(1-\exp\left(-\frac{\delta^2k}{2p}\right)\right).\\
\endalign
$$
One may chose $\delta=\sqrt{\frac{p}{2k}},$ which implies that
$$ E\leq
\left(\frac{20\delta}{\epsilon}\exp\left(-\frac{\epsilon^2k}{2p}\right)\right)I
\leq \left(\frac{p^4}{k^4}\right)I,$$ completing the proof.

\enddemo
\demo{ Proof of Proposition 3.2} If $-\pi/2\leq \omega\leq \pi/2$,
then
$$ \omega=\sin\omega + O\left(\sin^3\omega\right).$$  Also since
$\displaystyle{\cos\arg\left(1-\frac{e^{i\theta}}{p}\right)=\frac{1-\frac{\cos\theta}{p}}
{\left|1-\frac{e^{i\theta}}{p}\right|}}>0$, and
$\displaystyle{\sin\arg\left(1-\frac{e^{i\theta}}{p}\right)=\frac{-\sin
\theta}{p\left|1-\frac{e^{i\theta}}{p}\right|}}$, then
$$\arg\left(1-\frac{e^{i\theta}}{p}\right)=\frac{-\sin
\theta}{p\left|1-\frac{e^{i\theta}}{p}\right|}+O\left(\frac{\sin^3\theta}{p^3}\right).$$
This implies
$$ \cos\left(
2r\arg\left(1-\frac{e^{i\theta}}{p}\right)\right)= 1-
\frac{2r^2\sin^2\theta}{p^2\left|1-\frac{e^{i\theta}}{p}\right|^2}+
O\left(\frac{(r^4+r^2)\sin^4\theta}{p^4}\right).\tag{3.7}$$ Now
suppose that $p\leq k/\sqrt{\log k}$, in this case one has

$$\align
&\frac{1}{2\pi}\int_{-\pi}^{\pi}\left|1-\frac{e^{i\theta}}{p}\right|^{-2k}\cos\left(
2r\arg\left(1-\frac{e^{i\theta}}{p}\right)\right)d\theta=
\frac{1}{2\pi}\int_{-\pi}^{\pi}\left|1-\frac{e^{i\theta}}{p}\right|^{-2k}d\theta\\
&
-\frac{r^2}{p^2\pi}\int_{-\pi}^{\pi}\sin^2\theta\left|1-\frac{e^{i\theta}}{p}\right|^{-2k-2}d\theta+
O\left(\frac{r^4+r^2}{p^4}\int_{-\pi}^{\pi}\sin^4\theta\left|1-\frac{e^{i\theta}}{p}\right|^{-2k}d\theta\right).\\
\endalign
$$
Integrating by parts, we obtain
$$
\align
&\int_{-\pi}^{\pi}\sin^2\theta\left|1-\frac{e^{i\theta}}{p}\right|^{-2k-2}d\theta=
\int_{-\pi}^{\pi}\sin^2\theta\left(1-\frac{2\cos\theta}{p}+\frac{1}{p^2}\right)^{-k-1}d\theta\\
&=\frac{p}{2k}\int_{-\pi}^{\pi}\cos\theta\left(1-\frac{2\cos\theta}{p}+\frac{1}{p^2}\right)^{-k}d\theta
-\frac{p}{2k}\sin\theta\left(1-\frac{2\cos\theta}{p}+\frac{1}{p^2}\right)^{-k}\Bigg|_{-\pi}^{\pi}\\
&=\frac{p}{2k}\int_{-\pi}^{\pi}\cos\theta\left(1-\frac{2\cos\theta}{p}+\frac{1}{p^2}\right)^{-k}d\theta.\\
\endalign
$$
Further by Lemma 3.3, taking $\epsilon=
4\sqrt{\frac{p}{k}\log\left(\frac{k}{p}\right)}$ there,  we get
$$
\align
&\int_{-\pi}^{\pi}\cos\theta\left(1-\frac{2\cos\theta}{p}+\frac{1}{p^2}\right)^{-k}d\theta
=\left(1+O\left(\frac{p^4}{k^4}\right)\right)
\int_{-\epsilon}^{\epsilon}\cos\theta\left(1-\frac{2\cos\theta}{p}+\frac{1}{p^2}\right)^{-k}d\theta\\
&=\int_{-\epsilon}^{\epsilon}\left(1-\frac{2\cos\theta}{p}+\frac{1}{p^2}\right)^{-k}d\theta
\left(1+O\left(\epsilon^2+\frac{p^4}{k^4}\right)\right)\\
&=\int_{-\pi}^{\pi}\left(1-\frac{2\cos\theta}{p}+\frac{1}{p^2}\right)^{-k}d\theta
\left(1+O\left(\frac{p}{k}\log\left(\frac{k}{p}\right)\right)\right).\\
\endalign
$$
So we deduce that
$$\int_{-\pi}^{\pi}\sin^2\theta\left|1-\frac{e^{i\theta}}{p}\right|^{-2k-2}d\theta=
\frac{p}{2k}\int_{-\pi}^{\pi}\left|1-\frac{e^{i\theta}}{p}\right|^{-2k}d\theta
\left(1+O\left(\frac{p}{k}\log\left(\frac{k}{p}\right)\right)\right).$$
Moreover following the same ideas, and integrating by parts, we get
$$
\align
\int_{-\pi}^{\pi}\sin^4\theta\left|1-\frac{e^{i\theta}}{p}\right|^{-2k}d\theta&=
\frac{3p}{2(k+1)}\int_{-\pi}^{\pi}\sin^2\theta\cos\theta\left(1-\frac{2\cos\theta}{p}+\frac{1}{p^2}\right)^{-k+1}d\theta\\
&\ll
\frac{p^2}{k^2}\int_{-\pi}^{\pi}\left|1-\frac{e^{i\theta}}{p}\right|^{-2k}d\theta.\\
\endalign
$$
Thus by (3.7), the RHS of (3.6) equals
$$
\frac{1}{2\pi}\int_{-\pi}^{\pi}\left|1-\frac{e^{i\theta}}{p}\right|^{-2k}d\theta
\left(1-
\frac{r^2}{pk}+O\left(\frac{r^2}{k^2}\log\left(\frac{k}{p}\right)+\frac{r^4+r^2}{p^2k^2}
\right)\right).\tag{3.8}
$$
Now for the case $p\geq k/\sqrt{\log k}$, by (3.7) we use the
following estimate for the RHS of (3.6)
$$
\frac{1}{2\pi}\int_{-\pi}^{\pi}\left|1-\frac{e^{i\theta}}{p}\right|^{-2k}d\theta\left(1+
O\left(\frac{r^2}{p^2}\right)\right).\tag{3.9}$$ Finally upon using
(3.5), (3.6) and the estimates (3.8) and (3.9) for the appropriate
cases, we deduce that
$$
\align
&\sum_{n=1}^{\infty}\frac{d_{k-r}(n)d_{k+r}(n)}{n^2}=\sum_{n=1}^{\infty}\frac{d_k^2(n)}{n^2}
\prod_{p\geq k/\sqrt{\log
k}}\exp\left(O\left(\frac{r^2}{p^2}\right)\right)\\
& \prod_{p\leq k/\sqrt{\log k}}\exp\left(-
\frac{r^2}{pk}+O\left(\frac{r^2}{k^2}\log\left(\frac{k}{p}\right)+\frac{r^4+r^2}{p^2k^2}
\right)\right)\\
&= \exp\left(-r^2\frac{\log_2 k}{k}-\frac{c_0r^2}{k}
+O\left(\frac{r^2}{k\sqrt{\log
k}}+\frac{r^4}{k^2}\right)\right)\sum_{n=1}^{\infty}\frac{d_k^2(n)}{n^2},
\endalign
$$
completing the proof.
\enddemo
\head 4. Moments of $\zeta(1+it)$ \endhead

\noindent In this section we prove Theorem 2 together with a result
on moments of short Euler products. We begin by the proof of Theorem
2

\demo{ Proof of Theorem 2} First by Lemma 2.5, for $t$ large
enough we have
$$\zeta(1+it)^z=
\sum_{n=1}^{\infty}\frac{d_z(n)}{n^{1+it}}e^{-n/Z}
+O\left(\frac{1}{t}\right),$$ where $z$ is any complex number such
that $|z|\leq (\log t)^2$, and $Z=\exp((\log t)^{10})$.

\noindent Now let $x=\max\{|z_1|,|z_2|\}$, and $k=[x]+1$. Therefore
we have
$$
\align \frac{1}{T}\int_T^{2T}\zeta(1+it)^{z_1}&\zeta(1-it)^{z_2}dt\\
&= \sum_{m,n\geq
1}\frac{d_{z_1}(n)d_{z_2}(m)e^{-(m+n)/Z}}{mn}\frac{1}{T}\int_T^{2T}\left(\frac{m}{n}\right)^{it}dt
+ E_1,\tag{4.1}\\
\endalign
$$ where $$ E_1 \ll
\frac{1}{T}\sum_{n=1}^{\infty}\frac{d_k(n)}{n}e^{-n/Z}\ll
\frac{(\log 3Z)^{k}}{T} ,\tag{4.2}$$ by (3.2).
 The series in the RHS of (4.1) includes diagonal terms
$m=n$ which contribute as the main term,  and off-diagonal terms
$m\neq n$ which contribute as an error term, as we shall prove
later. The diagonal terms contribution  equals
$$ \sum_{n\geq
1}\frac{d_{z_1}(n)d_{z_2}(n)e^{-2n/Z}}{n^2}= \sum_{n\geq
1}\frac{d_{z_1}(n)d_{z_2}(n)}{n^2} +E_2,$$ where
$$ E_2\ll \frac{1}{\sqrt{Z}}\sum_{n\geq
1}\frac{d_{k}(n)^2}{n^{3/2}}\leq \frac{1}{\sqrt{Z}}\sum_{n\geq
1}\frac{d_{k^2}(n)}{n^{3/2}}=
\frac{\zeta(3/2)^{k^2}}{\sqrt{Z}},\tag{4.3}$$ knowing that
$1-e^{-t}\leq 2\sqrt{t}$ for all $t>0$.

\noindent For the off-diagonal terms, we divide the sum into four
parts:

\noindent a) $m,n\leq \sqrt T$, b) $m\geq n+n/\sqrt{T}$, c) $n\geq
m+m/\sqrt{T}$, and d) whatever remains. For the three first cases we
use the following inequality
$$\frac{1}{T}\int_T^{2T}\left(\frac{m}{n}\right)^{it}dt\ll
\frac{1}{T|\log(m/n)|} \leq T^{-1/2},\tag{4.4}$$ which holds since
$|\log(1-c)|=-\log(1-c)>c$ for any real number $0<c<1$.
 Thus by (3.2) the contribution of such terms is
  $$ E_3 \ll
  \frac{1}{T^{1/2}}\left(\sum_{n=1}^{\infty}\frac{d_k(n)}{n}e^{-n/Z}\right)^2\leq
\frac{(\log 3Z)^{2k}}{T^{1/2}}.\tag{4.5}$$ It remains then, to bound
the contribution from the last part $E_4$. we have
$$ |E_4|\leq \sum\Sb n+n/\sqrt{T}>m>n>\sqrt{T}\\ \hbox{or} \  m+m/\sqrt{T}>n>m>\sqrt{T}\endSb
 \frac{d_{|z_1|}(n)d_{|z_2|}(m)e^{-(n+m)/Z}}{mn}.$$ Let $y= \exp(\log
 T/\log_2 T)$. By Proposition 3.1 and (3.2), we get
 $$
 \align
 |E_4| &\leq 2
 \sum_{m>\sqrt{T}}\frac{d_k(m)}{m}e^{-m/Z}\sum_{m<n<m+m/\sqrt{T}}\frac{d_k(n)}{n}e^{-n/Z}\\
  & \leq 2
  \left(\sum_{m>\sqrt{T}}\frac{d_k(m)}{m}e^{-m/Z}\right)\left(\sup_{r>\sqrt{T}}
  \sum_{r<n<r+r/\sqrt{T}}\frac{d_k(n)}{n}e^{-n/Z}\right)\\
  &\leq 2(\log 3Z)^k\frac{(\log 3Z)^{k}}{y}=\frac{2(\log
  3Z)^{2k}}{y}.\\
 \endalign
 $$
This gives along with (4.2), (4.3) and (4.5), the following bound
$$E_1+E_2+E_3+E_4 \ll (2\log 3Z)^{2k}/y\leq y^{-1/2},$$  which achieves the proof.

\enddemo

Now to prove Theorem 1, we need a similar result to
Theorem 2, but for general short Euler products of degree 1. Indeed
we have

\proclaim{ Theorem 4.1} Let $T>0$ be large, $y\leq (\log T)^2$ a
real number, and $f$ a completely multiplicative function with
values on the unit circle $\Bbb U$ ($|f(n)|=1$ for all $n \in {\Bbb
N}$). Let
$$L_f(s,y):=\prod_{p\leq y}\left(1-\frac{f(p)}{p^s}\right)^{-1}=\sum_{n\in
S(y)}\frac{f(n)}{n^s},$$ where $S(y)=\{n\in {\Bbb N} : p|n\implies
p\leq y\}$. If $z_1,z_2$ are complex numbers verifying \noindent
$|z_1|,|z_2|\leq \displaystyle{\log T/(25\log_2 T\log_3 T)}$, then
$$ \frac{1}{T}\int_T^{2T}L_f(1+it,y)^{z_1}\overline{L_f(1+it,y)}^{z_2}dt =
\sum_{n\in S(y)} \frac{d_{z_1}(n)d_{z_2}(n)}{n^2} +
O\left(\exp\left(-\frac{\log T}{4\log_2 T}\right)\right).$$

\endproclaim
\demo{ Proof } Let $x=\max\{|z_1|,|z_2|\}$, and $k=[x]+1$. Then
$$ \frac{1}{T}\int_T^{2T}L_f(1+it,y)^{z_1}\overline{L_f(1+it,y)}^{z_2}dt =
\sum_{m,n\in
S(y)}\frac{d_{z_1}(n)d_{z_2}(m)f(n)\overline{f(m)}}{mn}\frac{1}{T}\int_T^{2T}\left(\frac{m}{n}\right)^{it}dt.
$$
\noindent In this series, the diagonal terms $m=n$ contribute
$$ \sum_{n\in S(y)}\frac{d_{z_1}(n)d_{z_2}(n)f(n)\overline{f(n)}}{n^2}=
\sum_{n\in S(y)}\frac{d_{z_1}(n)d_{z_2}(n)}{n^2}.$$

\noindent Furthermore we divide the off-diagonal terms into
two parts:

\noindent a) If $m,n\leq T^{3/4}$, and b) if $m>T^{3/4}$, or
$n>T^{3/4}$. Now for the first case we have
$$\frac{1}{T}\int_T^{2T}\left(\frac{m}{n}\right)^{it}dt\ll
\frac{1}{T|\log(m/n)|} \leq T^{-1/4},$$ since
$|\log(1-c)|=-\log(1-c)>c$ for any real number $0<c<1$. Thus the
contribution of such terms is bounded by
  $$
\frac{1}{T^{1/4}}\sum_{m,n\in
S(y)}\frac{d_{|z_1|}(n)d_{|z_2|}(m)}{mn}\leq
  \frac{1}{T^{1/4}}\left(\sum_{n\in S(y)}\frac{d_k(n)}{n}\right)^2
\leq \frac{(3\log y)^{2k}}{T^{1/4}}.\tag{4.6}
$$ Moreover the contribution from the second part is
bounded by
$$
 2\sum_{n,m\in S(y), n>T^{3/4}}
 \frac{d_{k}(n)d_{k}(m)}{mn}
 \leq
 2\left(T^{3/4}\right)^{-\alpha}\sum_{m\in S(y)}
 \frac{d_{k}(m)}{m}\sum_{n\in S(y)}
 \frac{d_{k}(n)}{n^{1-\alpha}},\tag{4.7}
$$ for all $\alpha>0$. We choose $\alpha=1/\log_2 T$.
Therefore
$$\sum_{n\in S(y)}
 \frac{d_{k}(n)}{n^{1-\alpha}}=\prod_{p\leq
 y}\left(1-\frac{1}{p^{1-\alpha}}\right)^{-k}=
 \exp\left(k\sum_{p\leq y}\frac{1}{p^{1-\alpha}}+O(k)\right)
 \leq \exp\left( 9 k\log_2 y\right).\tag{4.8}$$
Finally by  (4.6), (4.7) and (4.8) the contribution of the
off-diagonal part is at most
$$ \exp(-3\log T/4\log_2 T+ 10 k\log_2 y))\leq \exp(-\log T/4\log_2
T),$$ proving the Theorem.
\enddemo

 \head 5. Large values of $\zeta(1+it)$ in every direction\endhead

  In this section we prove Theorem 1.  For $s\in {\Bbb C}$, define $$L_{\psi}(s,y):=\prod_{p\leq
 y}\left(1-\frac{e^{-i\psi}}{p^{s}}\right)^{-1}.$$ We have
\proclaim{ Lemma 5.1} For $\theta\in(-\pi,\pi]$, and $y>0$ large
enough, let $\psi=\theta/\log_2 y$. Then for all  $t\in {\Bbb R}$,
we have
$$ L_{\psi}(1+it,y)=R_y\left(1+O\left(\frac{1}{\log_2 y}\right)\right)
 \Longleftrightarrow \zeta(1+it,y)= e^{i\theta}R_y\left(1+O\left(\frac{1}{\log_2 y}\right)\right).$$
\endproclaim

\demo{Proof} First we have
$$
\align \frac{1}{R_y}\prod_{p\leq
y}\left(1-\frac{e^{i\psi}}{p}\right)^{-1} &= \prod_{p\leq
y}\left(1-\frac{e^{i\psi}-1}{p-1}\right)^{-1}= \exp\left(-\sum_{p\leq
y}\log\left(1-\frac{e^{i\psi}-1}{p-1}\right)\right)\\
&=\exp\left(\sum_{p\leq y} \frac{i\psi}{p-1}+O\left(\psi^2\log_2
y\right)\right)=e^{i\theta}\left(1+O\left(\frac{1}{\log_2
y}\right)\right).\tag{5.1}
\endalign
$$
Now using that $(e^{i\psi})^m=1+O(m\psi)$ for all $m\in {\Bbb
N}$, we deduce that
$$
\align  \log\left(\prod_{p\leq
 y}\left(1-\frac{e^{-i\psi}}{p^{1+it}}\right)^{-1}\right)&-\log R_y=\sum_{p\leq y}\sum_{m\geq
1}\frac{(p^{-it}e^{-i\psi})^m-1}{p^mm}\\
&=\sum_{p\leq y}\sum_{m\geq
1}\frac{(p^{-it})^me^{-i\psi}-1}{p^mm}+O\left(\psi\sum_{p\leq
y}\sum_{m\geq 2}\frac{1}{p^m}\right)\\
&= e^{-i\psi}\sum_{p\leq y}\sum_{m\geq
1}\frac{(p^{-it})^m-e^{i\psi}}{p^mm}+O\left(\frac{1}{\log_2
y}\right).\tag{5.2}\\
\endalign
$$
Moreover by (5.1) we get
$$
\align \log \zeta(1+it,y)-\log\left(e^{i\theta}R_y\right)&=\log
\zeta(1+it,y)-\log\left(\prod_{p\leq
y}\left(1-\frac{e^{i\psi}}{p}\right)^{-1}\right)+O\left(\frac{1}{\log_2
y}\right)\\
&=\sum_{p\leq y}\sum_{m\geq
1}\frac{(p^{-it})^m-(e^{i\psi})^m}{p^mm}+O\left(\frac{1}{\log_2
y}\right)\\
&=\sum_{p\leq y}\sum_{m\geq 1}\frac{(p^{-it})^m-
 e^{i\psi}}{p^mm}+O\left(\psi\sum_{p\leq
y}\sum_{m\geq 2}\frac{1}{p^m}+\frac{1}{\log_2 y}\right)\\
&=\sum_{p\leq y}\sum_{m\geq 1}\frac{(p^{-it})^m-
 e^{i\psi}}{p^mm}+O\left(\frac{1}{\log_2 y}\right).\tag{5.3}\\
\endalign
$$
Finally the result follows upon taking absolute values of both (5.2)
and (5.3).
\enddemo

 Now we are ready to compute the moments (1.4). Indeed we
prove

 \proclaim{ Theorem 5.2 } Let $T>0$ be large, $y\leq \log T$ a real
 number, and $f$ a completely multiplicative function with values
on the unit circle $\Bbb U$. If $ k\leq y/(\log y)^2$ is a positive
integer and $\alpha=y/k$, then
$$
\align I(k):&=\frac{1}{T}\int_{T}^{2T}\left|L_{f}(1+it,y)
+R_y\right|^{2k} dt\\
&= (2R_y)^{2k}\exp\left(\frac{k}{\log k}\left(-\log(\alpha)+
O\left(1+\frac{(\log\alpha)^2}{\log
k}\right)\right)\right),\\
\endalign$$

\endproclaim
\demo{ Proof} First
$$
\align I(k)&= \frac{1}{T}\int_T^{2T}(L_{f}(1+it,y)
+R_y)^k(\overline{L_{f}(1+it,y)} +R_y)^kdt\\
& = \sum_{0\leq l,m\leq k}\binom kl\binom km
R_y^{2k-l-m}\frac{1}{T}\int_T^{2T}
L_{f}(1+it,y)^{l}\overline{L_{f}(1+it,y)}^{m}dt.\\
\endalign
$$ Therefore by
Theorem 4.1 we have
$$ I(k)= \sum_{0\leq l,m\leq k}\binom kl\binom km
R_y^{2k-l-m}\sum_{n\in S(y)}\frac{d_l(n)d_m(n)}{n^2} +
O\left((R_y+1)^{2k}\exp\left(-\frac{\log T}{4\log_2
T}\right)\right).$$ Now $k\leq y/(\log y)^2\leq \log T/(\log_2
T)^2$, and we know that $R_y\sim e^{\gamma}\log y$, thus
$$ (R_y+1)^{2k}\exp\left(-\frac{\log T}{4\log_2
T}\right)\leq \exp\left(-\frac{\log T}{5\log_2 T}\right).$$ We
divide the main term into two parts: central terms which correspond
to $l=m$, and non-central terms $l\neq m$.
\enddemo
\subhead { The lower bound }\endsubhead Since the contribution of
the non-central terms is positive, we have
$$ I(k)\geq \sum_{0\leq l\leq k}\binom kl^2
R_y^{2k-2l}\sum_{n\in S(y)}\frac{d_l(n)^2}{n^2} +
O\left(\exp\left(-\frac{\log T}{5\log_2 T}\right)\right).$$ Since
all the terms are positive, we consider only the contribution of
$l=[k/2]$. Then what remains only is to evaluate
 $\sum_{n\in S(y)}d_l(n)^2/n^2$. This has been done in [10] (Theorem
 3). Indeed Granville and Soundararajan proved that
$$\sum_{n\in S(y)}\frac{d_l(n)^2}{n^2}= \prod_{p\leq
l}\left(1-\frac{1}{p}\right)^{-2l}\exp\left(\frac{2l}{\log
l}\left(C+ O\left(\frac{l}{y}+\frac{1}{\log
l}\right)\right)\right),\tag{5.4}$$ where $C$ is the same constant
as (3). By (5.4) we get
$$ I(k)\geq \binom kl^2
R_y^{2k-2l}\prod_{p\leq
l}\left(1-\frac{1}{p}\right)^{-2l}\exp\left(\frac{2l}{\log
l}\left(C+ O\left(\frac{l}{y}+\frac{1}{\log
l}\right)\right)\right).$$ First by Stirling's formula we have
$\binom kl \gg 2^k/\sqrt{k}$. Moreover since $\log l/\log
y=1-\log(2\alpha)/\log y+O(1/\log^2 y)$, then
$$\prod_{p\leq
l}\left(1-\frac{1}{p}\right)^{-2l}=(R_y)^{2l}\exp\left(-2l\frac{\log(2\alpha)}{\log
y}\left(1+O\left(\frac{\log\alpha}{\log y}\right)\right)\right).
$$
Thus we deduce that
 $$ I(k)\geq \frac{1}{\sqrt{k}}(2R_y)^{2k}\exp\left(\frac{k}{\log
k}\left(-\log(\alpha)+ O\left(1+\frac{(\log\alpha)^2}{\log
k}\right)\right)\right),$$ which proves the lower bound.

\subhead{ The upper bound }\endsubhead
 \noindent Using Cauchy's inequality, we get that
$$
\align &\sum_{0\leq l, m\leq k}\binom kl\binom km
R_y^{2k-l-m}\sum_{n\in S(y)}\frac{d_l(n)d_m(n)}{n^2}\\
& \leq \sqrt{\left(\sum_{0\leq m\leq k}\sum_{0\leq l\leq k}\binom
kl^2
R_y^{2k-2l}\sum_{n\in S(y)}\frac{d_l(n)^2}{n^2}\right)^2}\\
&= (k+1) \sum_{0\leq l\leq k}\binom kl^2
R_y^{2k-2l}\sum_{n\in S(y)}\frac{d_l(n)^2}{n^2}.\\
\endalign
$$
Therefore $$ I(k)\leq (k+2) \sum_{0\leq l\leq k}\binom kl^2
R_y^{2k-2l}\sum_{n\in S(y)}\frac{d_l(n)^2}{n^2}.$$ Thus, by (5.4) we
deduce that
$$
\align I(k)&\leq (k+2) \sum_{0\leq l\leq k}\binom kl^2
R_y^{2k-2l}\prod_{p\leq
l}\left(1-\frac{1}{p}\right)^{-2l}\exp\left(\frac{2l}{\log
l}\left(C+ O\left(\frac{l}{y}+\frac{1}{\log
l}\right)\right)\right)\\
&\leq R_y^{2k}\sum_{0\leq l\leq k}\binom kl^2 \prod_{k\leq p\leq
y}\left(1-\frac{1}{p}\right)^{2l}\exp\left(O\left(\frac{k}{\log
k}\right)\right) \\
&= R_y^{2k}\sum_{0\leq l\leq k}\binom kl^2 \exp\left(\frac{2l}{\log
k}\left(-\log\alpha + O\left(\frac{(\log\alpha)^2}{\log
k}\right)\right)+O\left(\frac{k}{\log
k}\right)\right).\\
\endalign$$
Let $$f(l)=\binom kl^2\exp\left(-\frac{2l\log\alpha}{\log k}\right)
.$$ By Stirling's formula, one has
$$
\align \log f(l)=&2(k\log k-(k-l)\log(k-l)-l\log l)-\log 2\pi +\log
k-\log(k-l)-\log l \\
&-\frac{2l\log\alpha}{\log k}+o(1).\\
\endalign
$$ Differentiating the main term of this formula with respect to $l$, we deduce that the
maximum of $f$, occurs for
 $l=k/2(1+O(\log\alpha/\log k))$. Thus $$f(l)\leq
2^{2k}\exp\left(\frac{k}{\log k}\left(-\log\left(\alpha\right)+O
\left(1+\frac{(\log\alpha)^2}{\log k}\right)\right)\right),$$ which
implies the upper bound.

\demo{ Proof of Theorem 1} First, by Lemma 2.4 (with $A(y)=\log_2y$)
and Lemma 5.1, we have
$$
\align M(\theta,y)=&\text{meas}\{t\in [T,2T]:
L_{\psi}(1+it,y)=R_y\left(1+O\left(\frac{1}{\log_2
y}\right)\right)\} \\
&+ O\left(T\exp\left(-\frac{y\log_2y}{10^4\log y}\right)\right),\\
\endalign
$$ where $\psi=\theta/\log_2 y$ as in Lemma 5.1. Let $z$ be a
complex number verifying $|z|\leq 1$ and $|z+1|\geq 2-A\epsilon$,
for some positive constant $A$. Then $z=1+O(\sqrt{\epsilon})$ (where
the constant in the $O$ depends only on $A$). Moreover if
$z-1=O(\epsilon)$ then $|z+1|\geq 2-B\epsilon$ where $B$ depends
only on the constant in the $O$. Thus there exist some positive
constants $c_1$ and $c_2$ for which
$$ M_2+ O\left(T\exp\left(-\frac{y\log_2y}{10^4\log y}\right)\right)
\leq M(\theta,y) \leq M_1 + O\left(T\exp\left(-\frac{y\log_2
y}{10^4\log y}\right)\right),\tag{5.5}$$ where
$$ M_2:=\text{meas}\{t\in [T,2T] : |L_{\psi}(1+it,y)+R_y|\geq 2 R_y\left(1-\frac{c_2}{(\log_2
y)^2}\right)\},$$ and
$$ M_1:=\text{meas}\{t\in [T,2T] : |L_{\psi}(1+it,y)+R_y|\geq 2 R_y\left(1-\frac{c_1}{4\log_2
y}\right)\}.$$

\enddemo

\subhead The lower bound \endsubhead For a positive integer $k$, we
have
$$ TI(k)\leq (2R_y)^{2k}M_2+
(2R_y)^{2k}\exp\left(-2c_2\frac{k}{(\log_2 y)^2}\right)(T-M_2).$$
Now by Theorem 5.2, if $ k\leq y/(\log y)^2$, and $\alpha=y/k$, we
get
$$ T\exp\left(\frac{k}{\log
k}\left(-\log(\alpha)+ O\left(1+\frac{(\log\alpha)^2}{\log
k}\right)\right)\right)- T\exp\left(-2c_2\frac{k}{(\log_2
y)^2}\right) \leq M_2.$$ Choosing $k=[\exp(\log y- c_2\log y/(\log_2
y)^2)]$, we deduce
$$ M_2\geq 2T\exp\left(-\exp\left(\log y-c_2\frac{\log y}{(\log_2
y)^2}\right)\right).\tag{5.6}$$

\subhead The upper bound \endsubhead Similarly for a positive
integer $k$ we have
$$ (2R_y)^{2k}\exp\left(-c_1\frac{k}{\log_2 y}\right)M_1\leq
TI(k).$$ Then if $ k\leq y/(\log y)^2$, and $\alpha=y/k$, we get by
Theorem 5.2 that
$$ M_1\leq T\exp\left(c_1\frac{ k}{\log_2
y}+\frac{k}{\log
k}\left(-\log(\alpha)+O\left(1+\frac{(\log\alpha)^2}{\log k}\right)\right)\right).$$
Now by choosing $k=[\exp(\log y-2c_1\log y/\log_2 y)]$, we have that
$$ M_1\leq \frac{1}{2}T\exp\left(-\exp\left(\log y-c_1\frac{\log y}{\log_2
y}\right)\right).\tag{5.7}$$ Finally from (5.5), (5.6) and (5.7)
we deduce that
$$T\exp\left(-\exp\left(\log y-c_2\frac{\log y}{(\log_2
y)^2}\right)\right)\leq M(\theta,y)\leq T\exp\left(-\exp\left(\log
y-c_1\frac{\log y}{\log_2 y}\right)\right).$$

\head 6. Random Euler products and their distribution \endhead

\noindent We define $L(1,X,y):=\prod_{p\leq
y}\left(1-X(p)/p\right)^{-1}$. By the Central Limit Theorem,
$L(1,X,y)$ converges to $L(1,X)$ with probability $1$, as $y\to
\infty$. However we want a more accurate result which quantify the
rate of this convergence. Let $\Omega$ be the probability space on which
$\{X(p)\}_{p \ \text{prime}}$ are defined. For a real number $y>2$, define
$$ D(y):=\left\{\omega\in \Omega : L(1,X(\omega))=L(1,X(\omega),y)
\left(1+O\left(\frac{1}{\log y}\right)\right)\right\}.$$ Then we
prove

\proclaim{ Lemma 6.1} Let $y$ be large. We have
$$ 1-\text{Prob}(D(y))\ll\exp\left(-\frac{y}{e\log y}\right).$$
\endproclaim

\demo{ Proof} First we have
$$ \prod_{p>
y}\left(1-\frac{X(p)}{p}\right)^{-1} =\exp\left(\sum\Sb p>y\\ n\geq
1\endSb\frac{X(p)^n}{np^n}\right)=
\exp\left(\sum_{p>y}\frac{X(p)}{p}+O\left(\frac{1}{y}\right)\right).$$
Moreover
$$ {\Bbb E}\left(\left|\sum_{p>y}\frac{X(p)}{p}\right|^{2k}\right)={\Bbb E}\left(\sum\Sb p_1,...,p_k,q_1,...,q_k\\
p_i,q_j>y\endSb \frac{X(p_1)...X(p_k)\overline{X(q_1)...
X(q_k)}}{p_1...p_kq_1...q_k}\right).$$ Now if $p_1...p_k=q_1...q_k$
then ${\Bbb E}(X(p_1)...X(p_k)\overline{X(q_1)... X(q_k)})=1$,
otherwise this expectation is $0$. This gives $$ {\Bbb
E}\left(\left|\sum_{p>y}\frac{X(p)}{p}\right|^{2k}\right) \ll k!
\left(\sum_{p>y}\frac{1}{p^2}\right)^{k}\leq \left(\frac{k}{y\log
y}\right)^k.$$ Thus
$$ \left(\frac{1}{\log y}\right)^{2k}\text{Prob}\left(\left|\sum_{p>y}\frac{X(p)}{p}\right|>\frac{1}{\log y}\right)
\leq {\Bbb
E}\left(\left|\sum_{p>y}\frac{X(p)}{p}\right|^{2k}\right)\ll
\left(\frac{k}{y\log y}\right)^k.$$ Finally we choose $k=y/(e\log
y)$, which implies the result.

\enddemo

 To prove Theorem 3, we have to understand the correlation
between the norm and the argument of short Euler products of degree
$1$. For $y>2$ define
$$P_y:=\log\prod_{p\leq y}\left(1-\frac{1}{p}\right)^{-1}.$$
We have

\proclaim {Lemma 6.2} Let $\theta\ll 1$ and $\{x(p)\}_{p\leq y}$
 a sequence of complex numbers on the unit circle, such that $\arg\prod_{p\leq y}\left(1-x(p)/p\right)^{-1}=\theta.$ If
$$ \left|\prod_{p\leq
y}\left(1-\frac{x(p)}{p}\right)^{-1}\right|=\prod_{p\leq
y}\left(1-\frac{1}{p}\right)^{-1}\exp(L),$$ then
$$ L \leq -\frac{\theta^2}{2P_y}.$$
Moreover if $y$ is large, there exists
some real $\psi$ verifying

\noindent $\psi=\theta/P_y +O(\theta/P_y^2)$, and such that
$$\prod_{p\leq
y}\left(1-\frac{e^{i\psi}}{p}\right)^{-1}=\exp(L+i\theta)\prod_{p\leq
y}\left(1-\frac{1}{p}\right)^{-1},$$ with
$$L=-\frac{\theta^2}{2P_y}+O\left(\frac{\theta^2}{P_y^2}\right).$$
\endproclaim

\demo{Proof} The first statement of the Lemma follows upon noting
that $\log\prod_{p\leq
y}\left(1-\frac{x(p)}{p}\right)^{-1}=P_y+L+i\theta$, and
$$ \left|\log\prod_{p\leq
y}\left(1-\frac{x(p)}{p}\right)^{-1}=\sum\Sb p\leq y\\k\geq 1\endSb
\frac{x(p)^k}{kp^k}\right|\leq P_y.$$  For the second statement, we
search for $\psi$ such that
$$\theta=\sum\Sb p\leq y\\k\geq 1\endSb \frac{\sin(k\psi)}{kp^k}.$$
By the uniform convergence of the last series, $\psi$ exists and we
have
$$\theta= P_y\sin\psi + \sum\Sb p\leq y\\k\geq 2\endSb
\frac{\sin(k\psi)-\sin\psi}{kp^k}=\psi P_y+O(\psi+\psi^3P_y).$$ Thus
$ \psi=\theta/P_y +O(\theta/P_y^2),$ and finally
$$ L= \sum\Sb p\leq y\\k\geq 1\endSb \frac{\cos(k\psi)-1}{kp^k}
= (\cos\psi-1)P_y +\sum\Sb p\leq y\\k\geq 2\endSb
\frac{\cos(k\psi)-\cos\psi}{kp^k}=-\frac{\psi^2}{2}P_y+O(\psi^2),$$
which completes the proof.
\enddemo

\demo{ Proof of Theorem 3}  For $c>0$, and $y$ large enough we
define the following sets

$$
 B_+(c,\tau,y,\theta)= \left\{X\in  \Omega : |L(1,X,y)|>e^{\gamma}\tau
\left(1+\frac{c}{\log y}\right)\text{and} \ |\arg
L(1,X,y)|>\theta+\frac{c}{\log y}\right\},
$$

$$
 B_-(c,\tau,y,\theta)= \left\{X\in  \Omega : |L(1,X,y)|>e^{\gamma}\tau
\left(1-\frac{c}{\log y}\right)\text{and} \ |\arg
L(1,X,y)|>\theta-\frac{c}{\log y}\right\}.
$$
\enddemo
\subhead {The upper bound }\endsubhead If $c$ is a sufficiently
large constant, we get
$$ \Phi(\tau,\theta)\leq  \text{Prob}(B_-(c,\tau,y,\theta)\cap D(y)) + \text{Prob}(
D^c(y)).\tag{6.1}$$ Let $C_3>0$ be a suitably large constant and
choose
$$ \tau=\log
y\exp\left(-\frac{\theta^2}{2P_y}+2C_3\frac{\theta^2}{P_y^2}\right).$$

\noindent Take $X\in B_-(C_3,\tau,y,\theta)$, and put
$L(1,X,y)=\prod_{p\leq
y}\left(1-\frac{1}{p}\right)^{-1}\exp(L+i\phi)$. Then

\noindent $|\phi|>\theta-\displaystyle{\frac{C_3}{\log y}}$, and
$$ |L(1,X,y)|> e^{\gamma}\tau\left(1-\frac{C_3}{\log y}\right)
> \prod_{p\leq y}\left(1-\frac{1}{p}\right)^{-1}\exp
\left(-\frac{\theta^2}{2P_y}+C_3\frac{\theta^2}{P_y^2}\right).
$$
Therefore $$ L>-\frac{\theta^2}{2P_y}+C_3\frac{\theta^2}{P_y^2}\geq
\frac{-\phi^2}{2P_y}+C_3\frac{\theta^2}{P_y^2}-C_3\frac{\theta}{P_y\log
y} > \frac{-\phi^2}{2P_y}.$$ This contradicts  Lemma 6.2, which
implies that $B_-(C_3,\tau,y,\theta)=\emptyset$. Thus from (6.1) and
Lemma 6.1 we deduce that
$$ \Phi(\tau,\theta)\leq \text{Prob}(D^c(y))\ll
\exp\left(-\frac{y}{e\log y}\right).$$ And finally replacing $y$ by
$\tau$, we get
$$ \Phi(\tau,\theta) \leq \exp\left(-\frac{\displaystyle{e^{\tau+\frac{\theta^2\tau}{2\log\tau}
-3C_3\frac{\theta^2\tau}{\log^2\tau}}}}{\tau}\right),
$$ as desired.

\subhead {The lower bound } \endsubhead  By Lemma 6.1, if $c$ is a
sufficiently large constant, then
$$
\align \Phi(\tau,\theta)&\geq \text{Prob}(B_+(c,\tau,y,\theta)\cap
D(y)) \geq \text{Prob}(B_+(c,\tau,y,\theta)) +
\text{Prob}(D(y))-1\\
&\geq \text{Prob}(B_+(c,\tau,y,\theta)) -\exp\left(-\frac{y}{3\log
y}\right).\tag{6.2} \endalign$$ Now put $X(p)=e^{i\theta_p}$, where
the ${\theta_p}$ are independent random variables uniformly
distributed on $(-\pi,\pi)$. Let
$\widetilde{\theta}=\theta\left(1+1/P_y\right)$. By Lemma 6.2, there
exists $\psi$ verifying $$
 \psi=\frac{\widetilde{\theta}}{P_y}+O\left(\frac{\widetilde{\theta}}{P_y^2}\right),$$
 and such that $$\prod_{p\leq
y}\left(1-\frac{e^{i\psi}}{p}\right)^{-1}=\exp(L+i\widetilde{\theta})\prod_{p\leq
y}\left(1-\frac{1}{p}\right)^{-1},$$ where
$$ L=-\frac{\widetilde{\theta}^2}{2P_y}+O\left(\frac{\widetilde{\theta}^2}{P_y^2}\right)
=-\frac{\theta^2}{2P_y}+O\left(\frac{\theta^2}{P_y^2}\right).$$ We
choose $X\in \Omega$ such that
$$
\left\{\aligned &\psi-\frac{1}{(\log y)^{3/2}}<\theta_p<
\psi+\frac{1}{(\log y)^{3/2}} \ \ \text{for } \ p\leq z, \\ & \text{
and } \ \theta_p \in (-\pi,\pi), \ \text {for} \ z<p\leq
y,\\
\endaligned \right. \tag{6.3}$$
where $ z=\displaystyle{\frac{y}{8\log_2y}}$.  In this case
$$\prod_{p\leq
y}\left(1-\frac{X(p)}{p}\right)^{-1}=\prod_{p\leq
y}\left(1-\frac{e^{i\psi}}{p}\right)^{-1}\exp\left(O\left(\frac{1}{(\log
y)^{3/2}}\sum_{p\leq z}\frac{1}{p}+\sum_{z<p\leq
y}\frac{1}{p}\right)\right).$$  And since
 $$\sum_{z<p\leq y}\frac{1}{p}\sim \log\left(\frac{\log y}{\log
z}\right)=O\left(\frac{\log_3 y}{\log y}\right),$$ then
$$\prod_{p\leq
y}\left(1-\frac{X(p)}{p}\right)^{-1}=\exp\left(
-\frac{\theta^2}{2P_y}+i\left(\theta+\frac{\theta}{P_y}\right)
+O\left(\frac{\theta^2}{P_y^2}+\frac{\log_3 y}{\log
y}\right)\right)\prod_{p\leq y}\left(1-\frac{1}{p}\right)^{-1}.$$
Let $C_4>0$ be a suitably large constant, and choose
$$ \tau=\log
y\exp\left(-\frac{\theta^2}{2P_y}-2C_4\frac{\theta^2}{P_y^2}\right).$$
In this case we have
$$ \left|\prod_{p\leq
y}\left(1-\frac{X(p)}{p}\right)^{-1}\right|>e^{\gamma}\tau\left(1+\frac{C_4}{\log
y}\right),$$ and
$$ \left|\arg\prod_{p\leq
y}\left(1-\frac{X(p)}{p}\right)^{-1}\right|>\theta+\frac{C_4}{\log
y}.$$ Thus considering only these $X$ which satisfy (6.3) we deduce
that
$$ \text{Prob}(B_+(c,\tau,y,\theta))\geq \left(\frac{2}{2\pi(\log
y)^{3/2}}\right)^{\pi(z)}\geq \exp\left(-\frac{y}{4\log
y}+O\left(\frac{y}{\log y^2}\right)\right).$$ Finally by (6.2) we
get
$$\Phi(\tau,\theta)\geq \exp\left(-\frac{y}{3\log
y}\right)\geq\exp\left(-\frac{\displaystyle{e^{\tau+\frac{\theta^2\tau}{2\log\tau}
+3C_4\frac{\theta^2\tau}{\log^2\tau}}}}{\tau}\right).$$ Thus upon
taking $c_3=3C_3$, and $c_4=3C_4$, we deduce the result.

\demo{ Proof of Theorem 5} For the upper bound, the proof is the
same as for Theorem 3, replacing Lemma 6.1 by Lemma 2.4 (taking
$A(y)=\log y$). For the lower bound we use Theorems 4A and 4B to
make (6.3) holds in the appropriate ranges, and follow the same
lines as with Theorem 3.
\enddemo

\head 7. Fourier analysis on the $n$-dimensional torus \endhead

\noindent  We begin by presenting the following construction due to
Barton-Montgomery-Vaaler [1]:

 \noindent Let $N\in
{\Bbb N}$. If $u,v$ are real numbers with $0<u<v<1$, we define the
modified characteristic function $\phi_{u,v}: {\Bbb R}/{\Bbb
Z}\rightarrow {\Bbb R}$ by

$$ \phi_{u,v}(x)=\left\{\aligned 1 &  \ \text{ if } \  u<x-n<v \text{
for
} \ n\in {\Bbb Z}, \\
\frac{1}{2}& \ \text{ if} \ u-x\in {\Bbb Z} \ \text{or} \ v-x\in
{\Bbb Z},\\
0& \ \text{ otherwise.}\endaligned \right.$$ Put
$\bold{u}=(u_1,u_2,...,u_N)$ and $\bold{v}=(v_1,v_2,...,v_N)$ where
$0<u_n<v_n<1$. If $\bold L=(L_1,...,L_N)\in {\Bbb N}^N$, we let
${\Cal B}(\bold L)$ to be the set of all functions

\noindent $\bold \Phi_{\bold u,\bold v}$ : $({\Bbb R}/{\Bbb
Z})^N\rightarrow {\Bbb R}$ of the form
$$\bold \Phi_{\bold u,\bold
v}(\bold x)=\prod_{n=1}^N\phi_{u_n,v_n}(x_n),$$ and such that
$(v_n-u_n)(L_n+1)\in {\Bbb N}$ for all $1\leq n\leq N$. The principal result of Barton-Montgomery-Vaaler is the
following:

\proclaim {Theorem 7.1} Let $\bold L=(L_1,...,L_N)\in {\Bbb N}^N$
and $\bold \Phi_{\bold u,\bold v}\in{\Cal B}(\bold L)$.

\noindent There exist trigonometric polynomials $A(\bold x)$,
$B(\bold x)$ and $C(\bold x)$ of $N$ variables, with Fourier
coefficients  supported on the lattice
$$\Cal L:=\Cal L(L)= \{ l\in {\Bbb Z}^N: |l_n|\leq L_n, \
n=1,2,...,N\},$$  such that
$$ \aligned & \hat C(0)
=\prod_{n=1}^N\left(1+\frac{1}{(v_n-u_n)(L_n+1)}\right)\prod_{n=1}^N(v_n-u_n),\\
&\hat A(0)
=\prod_{n=1}^N(v_n-u_n),\\
&\hat B(0)
=\left(\sum_{n=1}^N\frac{1}{(v_n-u_n)(L_n+1)}\right)\prod_{n=1}^N(v_n-u_n),\\
\endaligned
$$
and
$$ A(\bold x)-B(\bold x)\leq \bold \Phi_{\bold u,\bold v}(\bold
x)\leq C(\bold x) \qquad \text{for all } \ \bold x \in ({\Bbb
R}/{\Bbb Z})^N.$$
\endproclaim
  Our goal is to prove Theorem 4A in the best possible
uniform region for $N=\pi(y)$. To this end we prove the following
Lemma which establishes the optimal choice of the lattice $\Cal L$
and thus of the degrees of the trigonometric polynomials we use
later in the proof of Theorem 4A.

\proclaim{ Lemma  7.2}  If $N=o\left(\sqrt{\log T/\log_2
T}\right)$, as $T\to \infty$, and $\{\delta_n\}_{1\leq n\leq N}$ are
real numbers between  $0$ and $1$ such that $$\min_{1\leq n\leq
N}\delta_N>\delta:=\left(\frac{N}{\sqrt{\log T/\log_2
T}}\right)^{2/3},$$ then there exist positive integers
$L_1,L_2,...,L_N$ verifying
$$ p_1^{L_1}p_2^{L_2}...p_N^{L_N}\leq T^{1/2},\tag{7.1}$$
and $$\sum_{n=1}^N\frac{1}{\delta_n(L_n+1)}=o(1).\tag{7.2}
$$
Moreover if (7.1) and (7.2) hold for some positive integers
$L_1,L_2,...,L_N$, and any real numbers $\{\delta_n\}_{1\leq n\leq
N}$ between $0$ and $1$, then $N=o\left(\sqrt{\log T/\log_2
T}\right)$.
\endproclaim
\demo{Proof} Let $L=\left[\log T/2N\right]$. If  $L_i=\left[L/\log
p_i\right]$, then
$$ \sum_{n=1}^NL_n\log p_n\leq LN\leq \log T/2,$$ which implies (7.1). Moreover
$$\sum_{n=1}^N\frac{1}{\delta_n(L_n+1)}\ll \frac{1}{\delta
L}\sum_{n=1}^N\log p_n\ll \frac{N\log N}{\delta L}\ll \frac{N^2\log
N}{\delta \log T}\ll \frac{N^{4/3}\log N}{(\log T)^{2/3}(\log_2
T)^{1/3}}=o(1),$$ and so (7.2) holds.

\noindent Now suppose that there exist positive integers
$L_1,L_2,...,L_N$, and real numbers $\{\delta_n\}_{1\leq n\leq N}$
between  $0$ and $1$, which verify (7.1) and (7.2). Then
$$ \sum_{n=1}^N \frac{1}{L_n}=o(1), \quad \text{and} \quad
\sum_{n=1}^N L_n\log p_n \leq \frac{\log T}{2}.$$ Thus by Cauchy's
inequality, we have
$$ \left(\sum_{n=1}^N\sqrt{\log p_n}\right)^2\leq\left(\sum_{n=1}^N
\frac{1}{L_n}\right)\left(\sum_{n=1}^N L_n\log p_n\right)=o(\log
T).$$ Finally by partial summation we get
$$ N\sqrt{\log N}\ll\sum_{n=1}^N\sqrt{\log p_n}=o(\sqrt{\log
T}),$$ which implies the result.
\enddemo

To prove Theorem 4B we need the following Lemma

\proclaim{ Lemma 7.3} Assume Conjecture 2. Let $N\leq \log
T/(10\log_2 T)$, as $T\to\infty$. Put $L=[N(\log T)^2]$. If
$|l_i|\leq L$, where $\{l_i\}_{1\leq i\leq N}$ are integers not all
zero, then
$$ |l_1\log p_1+l_2\log p_2+...+l_N\log p_N|\geq T^{-1/2}.$$

\endproclaim
\demo{ Proof} Let $\epsilon=1/100$. Since $\{\log p\}_{p \ \text{prime}}$ are linearly independent over $\Bbb Q$, there exists a constant
$c>0$ such that
$$
\align |l_1\log p_1+l_2\log p_2+...+l_N\log p_N|&>
\frac{c^NL}{(L^Np_1...p_N)^{1+\epsilon}}\\
&>\exp(-(1+2\epsilon)N\log N -(1+\epsilon)N\log L)\\
&>\exp\left(-\frac{\log T}{2}\right).\\
\endalign
$$
\enddemo

 \demo{ Proof of Theorem 4A} Let $a_j<b_j$ be the endpoints of  $I_j$,  and $L_j$ be positive integers
 satisfying the conditions of Lemma 7.2. There exist integers $0\leq r_i,s_i\leq L_i+1$ such that
$u_j:=r_j/(L_j+1)\leq a_j\leq x_j:=(r_j+1)/(L_j+1)$ and
$y_j:=s_j/(L_j+1)\leq b_j\leq v_j:=(s_j+1)/(L_j+1)$. Thus for all
$(z_1,z_2,...,z_N)\in ({\Bbb R}/{\Bbb Z})^N$, we have
$$ \bold \Phi_{\bold x,\bold y}(\bold
z):=\prod_{j=1}^N\phi_{x_j,y_j}(z_j)\leq
\prod_{j=1}^N\phi_{a_j,b_j}(z_j)\leq \bold \Phi_{\bold u,\bold
v}(\bold z):=\prod_{j=1}^N\phi_{u_j,v_j}(z_j).$$ Moreover $\bold
\Phi_{\bold x,\bold y},\bold \Phi_{\bold u,\bold v}\in \Cal B(\bold
L)$. Hence
$$\int_T^{2T}\prod_{j=1}^N\phi_{x_j,y_j}\left(\left\{\frac{t\log p_j}{2\pi}\right\}\right)dt
\leq
M\leq\int_T^{2T}\prod_{j=1}^N\phi_{u_j,v_j}\left(\left\{\frac{t\log
p_j}{2\pi}\right\}\right)dt. $$ Let $C(\bold z)$ be the
trigonometric polynomial as in Theorem 7.1, which corresponds to
$\bold \Phi_{\bold u,\bold v}$. Thus
$$
\align M&\leq \int_T^{2T}C\left(\left\{\frac{t\log
p_1}{2\pi}\right\},\left\{\frac{t\log
p_2}{2\pi}\right\},...,\left\{\frac{t\log
p_N}{2\pi}\right\}\right)dt\\
& =\int_T^{2T}\sum_{l\in \Cal L}\hat C (l)\exp\big(it(l_1\log
p_1+...+\l_N\log p_N)\big)dt\\
& =\sum_{l\in \Cal L}\hat C
(l)\int_T^{2T}\exp\left(it\log\left(p_1^{l_1}...
p_N^{l_N}\right)\right)dt.
\endalign
$$
The diagonal term which corresponds to $l=0$, equals $T\hat
C(0)$. Since $L_1,...,L_N$ verify the assertion
(7.1) of Lemma 7.2, it follows that the off-diagonal terms
contribute at most
$$\sum_{0\neq l\in \Cal L}|\hat C
(l)| \frac{2}{\left|\log\left(p_1^{l_1}...
p_N^{l_N}\right)\right|}\leq \left(\prod_{n=1}^N3L_n\right)\hat
C(0)p_1^{L_1}...p_n^{L_n}\leq \hat
C(0)\left(p_1^{L_1}...p_n^{L_n}\right)^{3/2}\leq T^{3/4}\hat C(0)
.$$ Finally since the assertion
(7.2) holds for our choices of $\delta_j$, we have
$$
\align M&\leq T\hat
C(0)\left(1+O\left(T^{-\frac{1}{4}}\right)\right)\\
&=T\prod_{n=1}^N(v_n-u_n)\prod_{n=1}^N\left(1+\frac{1}{(v_n-u_n)(L_n+1)}\right)
\left(1+O\left(T^{-\frac{1}{4}}\right)\right)\\
&=T\left(\prod_{n=1}^N\delta_n\right)\exp\left(O\left(\sum_{n=1}^N\frac{1}{L_j}+\sum_{n=1}^N\frac{1}{\delta_jL_j}
\right)\right)\left(1+O\left(T^{-\frac{1}{4}}\right)\right)\\
&=T\left(\prod_{n=1}^N\delta_n\right)(1+o(1)).
\endalign
$$
For the lower bound, we follow the same lines using the
corresponding trigonometric polynomials $A(\bold z)$ and $B(\bold
z)$ for $\bold \Phi_{\bold x,\bold y}$, as in Theorem 7.1. Indeed we
have
$$
\align M&\geq \int_T^{2T}(A-B)\left(\left\{\frac{t\log
p_1}{2\pi}\right\},\left\{\frac{t\log
p_2}{2\pi}\right\},...,\left\{\frac{t\log
p_N}{2\pi}\right\}\right)dt\\
& = T( \hat A(0)-\hat
B(0))\left(1+O\left(T^{-\frac{1}{4}}\right)\right)\\
&=T\prod_{n=1}^N(y_n-x_n)\left(1-\sum_{n=1}^N\frac{1}{(y_n-x_n)(L_n+1)}\right)
\left(1+O\left(T^{-\frac{1}{4}}\right)\right)\\
&=T\left(\prod_{n=1}^N\delta_n\right)(1+o(1)).
\endalign
$$
This completes the proof.

\enddemo

\demo{ Proof of Theorem 4B} The proof is exactly the same as Theorem 4A, taking  $L_j=[N(\log
T)^2]$ and using Lemma 7.3 instead of Lemma 7.2.
\enddemo
\head 8. The normal distribution of $\arg\zeta(1+it)$
\endhead

\noindent First we prove the following Lemma which shows that the
dominant contribution to the $2k$-th moment of $|\zeta(1+it)|$ comes
from the values of $t$ for which $|\zeta(1+it)|\approx
e^{\gamma}\tau$, provided that $k=e^{\tau-1-C}$, where $C$ is
defined by (3).

\proclaim{ Lemma 8.1} Let $T$, $\tau$, $\epsilon$, $k$, and
$\Omega_T(\tau)$ be as in Theorem 6. We have

$$ \frac{1}{T}\int_{T}^{2T}|\zeta(1+it)|^{2k}dt=
\frac{1}{T}\int_{\Omega_T(\tau)}|\zeta(1+it)|^{2k}dt
\left(1+O\left(\exp\left(-\frac{2k}{(\log
k)^{3/2}}\right)\right)\right).$$
\endproclaim
\demo{Proof} Upon integrating by parts, we get
$$
\align\frac{1}{T}\int_{\{t\in [T,2T]: \
|\zeta(1+it)|<e^{\gamma}(\tau-\epsilon)\}} &|\zeta(1+it)|^{2k}dt
=-e^{2k\gamma}\int_0^{\tau-\epsilon} x^{2k}d\Phi_T(x)\\
&=
e^{2k\gamma}\left(-(\tau-\epsilon)^{2k}\Phi_T(\tau-\epsilon)+2k\int_0^{\tau-\epsilon}\Phi_T(x)
x^{2k-1}dx\right).\tag{8.1}\\
\endalign
$$
Similarly one has
$$
\align
\frac{1}{T}\int_{\{t\in [T,2T]: \ |\zeta(1+it)|>e^{\gamma}(\tau+\epsilon)\}}&|\zeta(1+it)|^{2k}dt\\
&=
e^{2k\gamma}\left((\tau+\epsilon)^{2k}\Phi_T(\tau+\epsilon)+2k\int_{\tau+\epsilon}^{\infty}\Phi_T(x)
x^{2k-1}dx\right),\tag{8.2}\\
\endalign
$$ and
$$ \frac{1}{T}\int_{T}^{2T}|\zeta(1+it)|^{2k}dt=
e^{2k\gamma}\left(2k\int_{0}^{\infty}\Phi_T(x)
x^{2k-1}dx\right).\tag{8.3}$$ In [10], Granville and Soundararajan
proved that
$$2k\int_{0}^{\infty}\Phi_T(x) x^{2k-1}dx= (\log
k)^{2k}\exp\left(\frac{2k}{\log k}\left(C+O\left(\frac{1}{\log
k}\right)\right)\right),$$ together with
$$ \int_0^{\tau-\epsilon}\Phi_T(x)
x^{2k-1}dx \ll \exp\left(-\frac{2k}{(\log
k)^{3/2}}\right)\int_{0}^{\infty}\Phi_T(x)
x^{2k-1}dx,\tag{8.4}$$ and
$$ \int_{\tau+\epsilon}^{\infty}\Phi_T(x)
x^{2k-1}dx\ll \exp\left(-\frac{2k}{(\log
k)^{3/2}}\right)\int_{0}^{\infty}\Phi_T(x)
x^{2k-1}dx.\tag{8.5}$$
By (2) we deduce that
$$
\align \frac{(\tau+\epsilon)^{2k}\Phi_T(\tau+\epsilon)}
{\displaystyle{2k\int_{0}^{\infty}\Phi_T(x) x^{2k-1}dx}}&=
\left(\frac{\tau-C-1}{\tau+\epsilon}\right)^{-2k}\exp\left(-\frac{2e^{\tau+\epsilon-C-1}}{\tau+\epsilon}-\frac{2kC}{\log
k}+ O\left(\frac{k}{(\log k)^{3/2}}\right)\right)\\
&=
\exp\left(\frac{2k(1+C+\epsilon)}{\tau+\epsilon}-\frac{2ke^{\epsilon}}{\tau+\epsilon}-\frac{2kC}{\log
k}+O\left(\frac{k}{(\log k)^{3/2}}\right)\right)\\
&=
\exp\left(\frac{2k}{\tau}\left(1+\epsilon-e^{\epsilon}\right)+O\left(\frac{k}{(\log
k)^{3/2}}\right)\right)\\
&\ll \exp\left(-\frac{2k}{(\log k)^{3/2}}\right).\\
\endalign
$$
Similarly we get
$$\frac{(\tau-\epsilon)^{2k}\Phi_T(\tau-\epsilon)}
{\displaystyle{2k\int_{0}^{\infty}\Phi_T(x)
x^{2k-1}dx}}\ll\exp\left(-\frac{2k}{(\log k)^{3/2}}\right).$$
Finally using equations (8.1)-(8.5), we deduce the result.

\enddemo

\demo{ Proof of Theorem 6 } Let $x$ be a fixed real number, and
define
$$ \Lambda'_T(k,x):=\{t\in[T,2T]: \
\displaystyle{\frac{\arg\zeta(1+it)}{\sqrt{\frac{\log_2
k}{2k}}}}<x\}.$$ We consider the following distribution function
$$\nu'_{T,k}(x):=\frac{\displaystyle{\int_{\Lambda'_T(k,x)}|\zeta(1+it)|^{2k}dt}}
{\displaystyle{\int_T^{2T}|\zeta(1+it)|^{2k}dt}}.$$ The
characteristic function of  $\nu'_{T,k}$ is
$$
\psi_{T,k}(\eta):=\frac{\displaystyle{\int_T^{2T}|\zeta(1+it)|^{2k}\exp\left(i\eta\frac{\arg\zeta(1+it)}{\sqrt{\frac{\log_2
k}{2k}}}\right)dt}}
{\displaystyle{\int_T^{2T}|\zeta(1+it)|^{2k}dt}}.$$ Let
$\xi=\displaystyle{\frac{\eta}{\sqrt{\frac{\log_2 k}{2k}}}}$. One
can see that
 $\exp(i\arg\zeta(1+it))=\zeta(1+it)^{1/2}/\zeta(1-it)^{1/2}$,
 which implies
 $$\psi_{T,k}(\eta)= \frac{\displaystyle{\frac{1}{T}\int_T^{2T}\zeta(1+it)^{k+\xi/2}\zeta(1-it)^{k-\xi/2}dt}}
{\displaystyle{\frac{1}{T}\int_T^{2T}|\zeta(1+it)|^{2k}dt}}.$$ Now uniformly for $|\xi|\leq k$, we have by Theorem 2
$$
\psi_{T,k}(\eta)=\frac{\displaystyle{\sum_{n=1}^{\infty}\frac{d_{k+\xi/2}(n)d_{k-\xi/2}(n)}{n^2}}}
{\displaystyle{\sum_{n=1}^{\infty}\frac{d_{k}^2(n)}{n^2}}}+O\left(\exp\left(-\frac{\log
T}{2\log_2 T}\right)\right).$$ Finally by Proposition 3.2, and
replacing $\xi$ by $\eta$, we deduce that uniformly for
$|\eta|\leq\sqrt{\frac{\log_2k}{2}}$, we have
$$\psi_{T,k}(\eta)= \exp\left(-\frac{\eta^2}{2}-\frac{c_0\eta^2}{2\log_2k}
+O\left(\frac{\eta^2}{\sqrt{\log k}}+\frac{\eta^4}{\log_2^2
k}\right)\right)+O\left(\exp\left(-\frac{\log T}{2\log_2
T}\right)\right).$$ Let $\nu(x)=
\frac{1}{\sqrt{2\pi}}\int_{-\infty}^{x}e^{-y^2/2}dy$ be the normal
distribution function, and $\psi(\eta)=e^{-\eta^2/2}$  its
characteristic function. Then by the Berry-Esseen Theorem (Berry
[2], Esseen [9]),
$$ |\nu'_{T,k}(x)-\nu(x)|\leq K\int_{-R}^{R}\frac{|\psi_{T,k}(\eta)-\psi(\eta)|}{\eta}d\eta+\frac{B}{R},$$
for all $R>0$,  where $B$ and $K$ are absolute constants. We take
$R=\sqrt{\frac{\log_2k}{2}}$, which implies that
$$|\nu'_{T,k}(x)-\nu(x)|\ll
\int_{-R}^{R}\frac{e^{-\eta^2/2}\eta^2}{\eta\log_2k}d\eta+\frac{1}{\sqrt{\log_2k}}\ll
\frac{1}{\sqrt{\log_2k}}.$$ Finally by Lemma 8.1, we have
$$ \frac{1}{T}\int_{\Omega_T(\tau)}|\zeta(1+it)|^{2k}dt=
\frac{1}{T}\int_T^{2T}|\zeta(1+it)|^{2k}dt
\left(1+O\left(\exp\left(-\frac{2k}{(\log
k)^{3/2}}\right)\right)\right),$$ and
$$
\align
\frac{1}{T}\int_{\Lambda'_T(k,x)}|\zeta(1+it)|^{2k}dt&-\frac{1}{T}\int_{\Lambda_T(\tau,x)}|\zeta(1+it)|^{2k}dt\\
&\ll \exp\left(-\frac{2k}{(\log
k)^{3/2}}\right)\frac{1}{T}\int_T^{2T}|\zeta(1+it)|^{2k}dt\\
\endalign
$$
Therefore
$$\nu_{T,\tau}(x)=\nu'_{T,k}(x)+O\left(\exp\left(-\frac{2k}{(\log
k)^{3/2}}\right)\right),$$ which completes the proof.
\enddemo

\head 9. Analogous results for $L(1,\chi)$ \endhead

\noindent In this section we present the analogous results for
$L(1,\chi)$. Although we expect the behavior of the sets of values
of $\zeta(1+it)$ and these of $L(1,\chi)$ should be the same, one
should note that there are some differences between these two sets.
Indeed the first set is continuous and the moments are integrals,
while the second one is discrete and the moments are sums. Also an
extra difficulty in the case of $L(1,\chi)$, is the possible
existence of {\it Landau-Siegel} zeros, corresponding to exceptional
{\it Siegel} characters $\chi$ defined as follows

$$ \chi \  \hbox{mod} \  q: \hbox{ there exists } s \hbox{ with }
\ \text{Re} (s) \geq 1- \frac{c}{\log q (\hbox{Im} (s) +2)}  \hbox{ and } \
L(s,\chi)=0,$$ for some small constant $c>0$. Let $S$ be the set of
such characters. One expects this set to be empty, but what is known
unconditionally (see [5]), is that such characters are very rare.
Indeed each $\chi$ must be real (thus of order 2), and between any
two powers of $2$ there is at most one fundamental discriminant $D$
with $\left(\frac{D}{\cdot}\right)\in S$.  Throughout this section
$q$ will denote a large prime number. In this case there is at most
one exceptional character $\chi$ of conductor $q$.

Using similar ideas, we show the existence of large values of
$L(1,\chi)$ in every direction

\proclaim{ Theorem 9.1 } Fix $\theta\in (-\pi,\pi]$. If $1\ll y\leq
\log q/\log_2 q$ is a real number, let $N(\theta,y)$ be the number
of non-principal characters $\chi \notin S$ of conductor $q$ for
which
$$ L(1,\chi) =e^{i\theta}\prod_{p\leq
y}\left(1-\frac{1}{p}\right)^{-1}\left(1+O\left(\frac{1}{\log_2
y}\right)\right).$$ Then there exist two positive constants
$c_6,c_7$ (depending on the constant in the $O$) for which
$$\phi(q)\exp\left(-y^{1-c_6/(\log_2
y)^2}\right)\leq N(\theta,y)\leq \phi(q)\exp\left(-y^{1-c_7/\log_2
y}\right).$$
\endproclaim

\demo{ Proof} We follow exactly the proof of Theorem 1: first we
prove the analogue of Theorem 4.1 to get asymptotic for moments of
short Euler products $\prod_{p\leq y} (1-f(p)\chi(p)/p)^{-1}$, where
$f$ is a completely multiplicative function with values on the unit
circle. Then we prove the analogue of Lemma 5.1, replacing $p^{-it}$
by $\chi(p)$ (the proof is the same since $|\chi(p)|=1$). What
remains is to prove the analogue of Lemma 2.4, which can be done
using the zero free region and zero density estimates of
$L(s,\chi)$, if $\chi\notin S$.
\enddemo

 As mentioned in the introduction, using a different approach , Granville and
Soundararajan (unpublished) proved the existence of large values
(and small ones) in every direction. Indeed what they established is
the following \proclaim { Theorem A (Gr-S)} If  $z$ is any complex
number such that

$$ \frac{\pi^2}{6 e^{\gamma}\log_2 q}\left(1+O\left(\frac{1}{\log_3 q
}\right)\right)\leq |z|\leq e^{\gamma} \log_2 q
\left(1+O\left(\frac{1}{\log_3 q }\right)\right),$$ then the number of
non-principal characters $\chi \notin S$ of conductor $q$ for which
$$ L(1,\chi)= z \left( 1+O\left( \frac{\log_3q}{\log_2
q}\right)\right),$$ is at least $q^{1-1/\log_2 q}$.

\endproclaim

Also in their unpublished draft, they proved an analogue of Theorem 2 for
complex moments of $L(1,\chi)$

\proclaim { Theorem B (Gr-S)} Fix $\epsilon>0$ and suppose that $q$
is a sufficiently large integer. Let $H$ be a subgroup of the
character group $G$ for $({\Bbb Z}/q{\Bbb Z})^*$ with $[G:H]\ll
\exp\left(\log^{\epsilon/2}q\right)$. Assume that there is an
integer $r\leq \log^{1-\epsilon}q$ for which $\chi^r\in H$ for all
$\chi\in G$. If $z_1$ and $z_2$ are complex numbers with
$|z_1|,|z_2|\leq \displaystyle{\log q/r(\log_2 q)^3}$, and $\xi$
is any character in $G$ then, we have uniformly
$$ \frac{1}{|H_{\xi}|}\sum\Sb \chi\in H_{\xi}\endSb L(1,\chi)^{z_1}L(1,\overline{\chi})^{z_2} =
\sum\Sb n=1\\ (n,q)=1\endSb ^{\infty}
\frac{d_{z_1}(n)d_{z_2}(n)}{n^2} + o(1),$$ where $H_{\xi}$ is the
set of characters in $\xi H$ of order $>1$, not belonging to $S$.
\endproclaim
\noindent If we restrict our selves to the case of $q$ prime then using a similar approach as in Theorem 2, we
have

\proclaim { Theorem 9.2} Let $q$ be a large prime. Then uniformly
for all complex numbers $z_1,z_2$ in the region
 $|z_1|,|z_2|\leq \displaystyle{\log q/50(\log_2 q)^2}$,
we have
$$ \frac{1}{\phi(q)}\sum\Sb\chi \ (\text{mod} \ q)\\ \chi\neq 1, \chi\notin S\endSb L(1,\chi)^{z_1}L(1,\overline{\chi})^{z_2} =
\sum\Sb n=1\\ (n,q)=1\endSb ^{\infty}
\frac{d_{z_1}(n)d_{z_2}(n)}{n^2} + O\left(\exp\left(-\frac{\log
q}{2\log_2 q}\right)\right).$$
\endproclaim

\demo{ Proof} We follow the same lines as the proof of Theorem 2. First by Lemma 2.3 of [11] (analogue of Lemma 2.5 for
$\zeta(1+it)$), if $\chi$ is a non-principal character $(\hbox{mod }
q)$ not belonging to $S$, then
$$ L(1,\chi)^z=\sum_{n=1}^{\infty}\chi(n)\frac{d_z(n)}{n}e^{-n/Z}+O\left(\frac{1}{q}\right),$$
where $Z=\exp\left((\log q)^{10}\right)$ and $z$ is any complex
number with $|z|\leq (\log q)^2$. Let $k=\max\{[|z_1|]+1,
[|z_2|]+1\}$, we have then
$$
\align
 \frac{1}{\phi(q)}\sum\Sb\chi \ (\text{mod} \ q)\\ \chi\neq 1, \chi\notin S\endSb & L(1,\chi)^{z_1}L(1,\overline{\chi})^{z_2}\\
&= \sum_{n,m\geq
1}\frac{d_{z_1}(n)d_{z_2}(m)}{nm}e^{-(m+n)/Z}\frac{1}{\phi(q)}\sum\Sb\chi
\ (\text{mod} \ q)\\ \chi\neq 1, \chi\notin S\endSb
\chi(n)\overline{\chi(m)}+O\left(\frac{(\log 3 Z)^k}{q}\right).\tag{9.1}\\
\endalign
$$
We now extend the right side of (9.1) so as to include all
characters $(\hbox{mod } q)$. Since $q$ is prime, $S$ contains at
most one element, thus by (3.2) the contribution of characters of
$S$ together with the principal character is bounded by
$$ \frac{2}{\phi(q)}\left(\sum_{n\geq
1}\frac{d_k(n)}{n}e^{-n/Z}\right)^2\ll \frac{(\log 3 Z)^{2k}}{q}.$$
The
contribution from the diagonal terms $m=n$ is
$$ \sum\Sb n=1\\ (n,q)=1\endSb^{\infty}\frac{d_{z_1}(n)d_{z_2}(n)}{n^2}e^{-2n/Z}
= \sum\Sb n=1\\
(n,q)=1\endSb^{\infty}\frac{d_{z_1}(n)d_{z_2}(n)}{n^2} +
O\left(\frac{\zeta(3/2)^{k^2}}{\sqrt{Z}}\right),$$ by (4.3). Using the orthogonality relations for characters, we see that
the off-diagonal terms $m\neq n$ satisfy $m\equiv n \ (\hbox{mod }
q)$ and $(mn,q)=1$, which imply $\max (m,n)> q$. Thus the
contribution of these terms is bounded by
$$ 2\sum_{n=1}^{\infty}\frac{d_k(n)}{n}e^{-n/Z}\left(\max_{ b \
\text{mod} \ q}\sum\Sb m> q\\ m\equiv b \ \text{mod} \ q\endSb
\frac{d_k(m)}{m}e^{-m/Z}\right).\tag{9.2}$$ Now following the proof
of Proposition 3.1 (using induction on $k$), we can prove that
$$ \max_{ b \
\text{mod} \ q}\sum\Sb m> q\\ m\equiv b \ \text{mod} \ q\endSb
\frac{d_k(m)}{m}e^{-m/Z}\leq \frac{(\log 3 Z)^{k}}{y},\tag{9.3}$$
where $y=\exp(\log q/\log_2 q)$. Finally by (9.2) and (9.3) we
deduce the result.
\enddemo

 Using Fourier analysis on the $n$-dimensional torus, and
the construction of Barton-Montgomery-Vaaler [1], we proved the
uniform distribution of the values $\{ p^{it} : t\in [T,2T]\}_{p\leq
y}$. We can use exactly the same ideas to prove that the values $\{
\chi(p): \chi \text{ mod } q\}_{p\leq y}$ have the same behavior.
Indeed we have

\proclaim {Theorem 9.3} Let $2<y$ be a real number. For each $1\leq
j\leq \pi(y)$, let $I_j\subset (0,1)$ be an open interval of length
$\delta_j>0$. Define
$$ N(I_1,...,I_{\pi(y)})=N:=\left|\left\{ \chi \text { mod } q:  \left\{ \frac{\arg(\chi(p_j))}{2\pi} \right\}\in
I_j, \text{ for all } 1\leq j\leq \pi(y) \right\}\right|,$$ where
$p_j$ is the $j$-th smallest prime, and $\{\cdot\}$ denotes the
fractional part. We have
$$N\sim \phi(q) \prod_{j\leq \pi(y)}\delta_j,$$
uniformly for $y\leq \sqrt{\log q}/(\log_2 q)^2$, and
$\delta_j>(\log_2q)^{-5/3}$.
\endproclaim

\noindent One should note that there is no analogue of Theorem 4B
(where we assume Conjecture 2) in this case.

\demo{ Proof} The proof is exactly the same as Theorem 4A, noting
that
$$ \sum _{\chi \ (\text {mod } q)} A\left(
\frac{\arg(\chi(p_1))}{2\pi},...
,\frac{\arg(\chi(p_n))}{2\pi}\right)= \phi(q) \hat A(0),\tag{9.4}$$ if $A$ is
a trigonometric polynomial in $n$ variables, with Fourier
coefficients supported in a lattice
$$\Cal L= \{ l\in {\Bbb Z}^n: |l_i|\leq L_i,
i=1,2,...,n\},$$ with $p_1^{L_1}p_2^{L_2}...p_n^{L_n}\leq q$. This
follows from the orthogonality relation for characters and the fact
that
$$
\align &\sum _{\chi \ (\text {mod } q)} A\left(
\frac{\arg(\chi(p_1))}{2\pi},...
,\frac{\arg(\chi(p_n))}{2\pi}\right)\\
&=\sum _{\chi \ (\text {mod } q)} \sum_{l\in \Cal L} \hat A(l) \exp(
i(l_1 \arg(\chi(p_1))+...+l_n\arg(\chi(p_n))))= \sum_{l\in \Cal L}
\hat A(l)\sum _{\chi \ (\text {mod } q)}
\chi\left(\prod_{i=1}^np_i^{l_i}\right).\\
\endalign$$

\enddemo

 Finally we can use the same ideas in the proofs
of Theorems 5 and 6, to deduce analogous results for $L(1,\chi)$. Indeed define

$$\Phi_q(\tau):=\frac{1}{\phi(q)}|\{\chi \text{ mod } q, \chi\neq 1, \chi\notin S :
|L(1,\chi)|>e^{\gamma}\tau\}|, \ \text{ and }$$

$$\Phi_q(\tau,\theta):=\frac{1}{\phi(q)}|\{\chi \text{ mod } q, \chi\neq 1, \chi\notin S :
|L(1,\chi)|>e^{\gamma}\tau,
 \ |\arg L(1,\chi)|>\theta\}|.$$ In [10], Granville and Soundararajan proved that the asymptotic relation (2)
 holds also for $\phi_q(\tau)$. For $\Phi_q(\tau,\theta)$, similarly to Theorem 5 we prove

\proclaim{ Theorem 9.4} Let $q$ be a large prime number. There exist
two positive constants $c_8$ and $c_9$ such that
$$
\Phi_q(\tau,\theta)\leq
\exp\left(-\frac{\displaystyle{e^{\tau+\frac{\theta^2\tau}{2\log\tau}
-c_8\frac{\theta^2\tau}{\log^2\tau}}}}{\tau}\right),$$
uniformly for $1\ll\tau\leq \log_2 q$, and $(\log
\tau)\sqrt{\frac{\log_2 \tau}{\tau}}<\theta\ll 1$. And
$$
\Phi_q(\tau,\theta)\geq
\exp\left(-\frac{\displaystyle{e^{\tau+\frac{\theta^2\tau}{2\log\tau}
+c_9\frac{\theta^2\tau}{\log^2\tau}}}}{\tau}\right),$$
uniformly for $1\ll\tau\leq (\log_2 q)/2-2\log_3 q$ and
$(\log \tau)\sqrt{\frac{\log_2\tau}{\tau}}<\theta\ll 1$.
\endproclaim

\demo{ Proof } For the upper bound, the proof is the same as for
Theorem 5, replacing Lemma 6.1 by the analogue of
Lemma 2.4 for $L(1,\chi)$ (taking
$A(y)=\log y$). For the lower bound we use Theorem 9.3 to make (6.3)
holds and follow the same lines as with Theorem 5.

\enddemo

\noindent

\proclaim {Corollary 9.1} If $1\ll\tau\leq \log_2 q-\log_3 q$, then
for almost all characters $\chi \mod q$, with
$|L(1,\chi)|>e^{\gamma}\tau$, we have $|\arg L(1,\chi)|\leq (\log
\tau)\sqrt{\log_2\tau/\tau}$.
\endproclaim

 We prove also

\proclaim{ Theorem 9.5} Let $q$ be a large prime number,
$1\ll\tau\leq \log_2q-3\log_3q$ a real number,
$\epsilon=\tau^{-1/5}$ and $k=e^{\tau-1-C}$, where $C$ is defined by
(3). Let
$$ \Omega_q(\tau):=\{\chi \text{ mod } q, \chi\neq 1, \chi\notin S : \ e^{\gamma}(\tau-\epsilon)\leq
|L(1,\chi)|\leq e^{\gamma}(\tau +\epsilon)\},$$ and for a real
number $x$, let
$$ \Lambda_q(\tau,x):=\{\chi\in\Omega_q(\tau): \
\displaystyle{\frac{\arg
L(1,\chi)}{\sqrt{\frac{\log(\tau-1-C)}{2e^{\tau-1-C}}}}}<x\} \ \text{ and} \ \nu_{q,\tau}(x):=\frac{\displaystyle{\sum_{\chi\in\Lambda_q(\tau,x)}|L(1,\chi)|^{2k}}}
{\displaystyle{\sum_{\chi\in\Omega_q(\tau)}|L(1,\chi)|^{2k}}}.$$ Then  we have
$$\nu_{q,\tau}(x)=\frac{1}{\sqrt{2\pi}}\int_{-\infty}^{x}e^{-y^2/2}dy+O_x\left(\frac{1}{\sqrt{\log\tau}}\right).$$
\endproclaim

\demo{ Proof } The proof is exactly the same as Theorem 6, using
Theorem 9.2, along with Proposition 3.2 and the results of
Granville-Soundararajan [10] for the distribution of $|L(1,\chi)|$
(which are exactly the same as for $|\zeta(1+it)|$).
\enddemo


\Refs

\ref \key 1 \by J.T. Barton, H.L. Montgomery and J.D. Vaaler  \book
Note on a Diophantine inequality in several variables \publ Proc.
Amer. Math. Soc. \bf{129} \yr (2001), 337-345
\endref


\ref \key 2 \by A.C. Berry  \book The accuracy of the Gaussian
approximation to the sum of independent variables \publ Trans. Amer.
Math. Soc. \bf{49} \yr (1941), 122-136
\endref

\ref \key 3 \by E. Bombieri \book Le grand crible en th\'eorie
analytique des nombres \publ Ast\'erisque \bf{18} \yr (1987/1974),
103 pp
\endref

\ref \key 4 \by J. Cogdell and P. Michel \book On the complex
moments of symmetric power $L$-functions at $s=1$ \publ IMRN \bf{31}
\yr (2004), 1561-1617
\endref

\ref \key 5 \by H. Davenport \book Multiplicative number theory
\publ Springer Verlag, New York \yr 1980.
\endref

\ref \key 6 \by W. Duke \book Extreme values of Artin $L$-functions and class numbers
\publ Compositio Math. \bf{136} \yr (2003), no. 1, 103-115
\endref

\ref \key 7 \by P.D.T.A. Elliott \book On the size of $L(1,\chi)$
\publ J. reine angew. Math. \bf{236} \yr (1969), 26-36
\endref

\ref \key 8 \by P.D.T.A. Elliott \book On the distribution of the
values of quadratic $L$-series in the half-plane $\sigma>1/2$ \publ
Invent. Math. \bf{21} \yr (1973), 319-338
\endref

\ref \key 9 \by C.G. Esseen \book Fourier analysis of distribution
functions, A mathematical study of the Laplace-Gaussian law \publ
Acta Math. \bf{77} \yr (1945), 1-125
\endref


\ref \key 10 \by A. Granville and K. Soundararajan \book Extreme
values of $|\zeta(1+it)|$ \publ  The Riemann zeta function and related themes: papers in honour of Professor K. Ramachandra, Ramanujan Math. Soc. Lect. Notes Ser., {\bf 2} \yr (2006), 65-80
\endref

\ref \key 11 \by A. Granville and K. Soundararajan \book The
distribution of values of $L(1,\chi_d)$ \publ Geometric and Funct.
Anal \bf{13} \yr (2003),  992-1028
\endref

\ref \key 12 \by L. Habsieger  and E. Royer \book $L$-functions of automorphic forms and combinatorics: Dyck paths  \publ Ann. Inst. Fourier (Grenoble) \bf{54} \yr (2004), no.7,  2105-2141
\endref

\ref \key 13 \by S. Lang \book Elliptic curves: Diophantine analysis
\publ Grundlehren der Mathematischen Wissenschaften [Fundamental
Principles of Mathematical Sciences], {\bf 231}. Springer-Verlag,
Berlin-New York \yr 1978
\endref

\ref \key 14 \by Y.-K. Lau  and J. Wu \book Extreme values of symmetric power $L$-functions at $1$  \publ Acta Arith. \bf{126} \yr (2007), no.1,  57-76
\endref


\ref \key 15 \by J.E. Littlewood \book On the function $1/\zeta
(1+it)$ \publ Proc. London Math. Soc \bf{ 27} \yr (1928),  349-357
\endref


\ref \key 16 \by J.E. Littlewood \book On the class number of the
corpus $P(\sqrt{-k})$ \publ Proc. London Math. Soc \bf{27} \yr (1928), 358-372
\endref


\ref \key 17 \by J.Y. Liu, E. Royer and J. Wu \book On a conjeture
of Montgomery-Vaughan on extreme values of automorphic $L$-functions
at $1$. To appear at proceedings for the conference on "Anatomy of integers", Montreal \yr (2006)
\endref


\ref \key 18 \by H.L. Montgomery  \book Ten lectures in the
interface between Analytic Number Theory and Harmonic Analysis \publ
American Mathematical Society, Providence, RI \yr (1994)
\endref

\ref \key 19 \by H.L. Montgomery and R.C. Vaughan \book Extreme
values of Dirichlet L-functions at $1$ \publ "Number Theory in
Progress"(K. Gy\"ory, H. Iwaniec, J. Urbanowicz, eds.), de Gruyter,
Berlin  \yr (1999), 1039-1052
\endref

\ref \key 20 \by K.K. Norton  \book Upper bounds for sums of powers
of divisor functions \publ J. Number Theory
  \bf{40} \yr (1992),  60-85
\endref

\ref \key 21 \by E. Royer  \book Statistique de la variable al\'eaoire $L(sym^2 f,1)$ \publ Math. Ann.
  \bf{321} \yr (2001), no. 3,  667-687
\endref

\ref \key 22 \by E. Royer  \book Interpr\'etation combinatoire des moments n\'egatifs des valeurs de fonctions $L$ au bord de la bande critique  \publ Ann. Sci. \'Ecole Norm. Sup. (4)
  \bf{36} \yr (2003), no. 4,  601-620
\endref

\ref \key 23 \by E. Royer and J. Wu \book Taille des valeurs de fonctions $L$ de carr\'es sym\'etriques au bord de la bande critique \publ Rev. Mat. Iberoamericana
  \bf{21} \yr (2005), no. 1,  263-312
\endref

\ref \key 24 \by E. Royer and J. Wu \book Special values of symmetric power $L$-functions and Hecke eigenvalues  \publ J. Th\'eor. Nombres Bordeaux
  \bf{19} \yr (2007), no. 3,  703-753
\endref


\ref \key 25 \by E.C. Titchmarsh \book The theory of the Riemann
zeta-function \publ Oxford University Press, Oxford \yr 1986
\endref

\endRefs
\enddocument